\documentclass[pdftex,twoside,leqno]{article}

\usepackage[english]{babel}
\usepackage{graphicx}
\usepackage{hyperref}
\usepackage{bbold}
\usepackage{url}
\usepackage{amsmath}
\usepackage{amsfonts}
\usepackage{color}
\usepackage{verbatim}
\usepackage[utf8]{inputenc}

\usepackage{amsthm} 
\theoremstyle{plain}
\newtheorem{theorem}{Theorem}[section]
\newtheorem{lemma}{Lemma}[section]

\makeatletter
\newcommand{\toto@c@page}{}
\newcommand{\toto@addmarginpar}{}
\let\toto@addmarginpar\@addmarginpar
\renewcommand{\@addmarginpar}{
  \def\toto@c@page{\c@page}
  \c@page=1
  \toto@addmarginpar
  \c@page=\toto@c@page
}
\makeatother
\usepackage[textwidth=2.5cm,textsize=footnotesize]{todonotes}

\usepackage{listings}

\lstdefinelanguage{Coq}{
  comment=[n]{(*}{*)},
  otherkeywords={->, /\\, \\/, =>},
  morekeywords={Definition, End, Fixpoint, Goal, Lemma, Notation, Proof, Qed, Require, Section, Theorem, Time, Variable, apply, by, end, exists, forall, fun, intros, let, match, replace, ring, with},
  literate={=>}{{$\Rightarrow\,$}}1 {->}{{$\to\,$}}1 {/\\}{{$\land\,$}}1
  {<=}{{$\le\,$}}1 {>=}{{$\ge\,$}}1 {<>}{{$\neq\,$}}1
}

\lstdefinestyle{ACSL}{
  language=C,
  literate={==}{{$=\,$}}1 {<=}{{$\le\,$}}1 {>=}{{$\ge\,$}}1
  {==>}{{$\Rightarrow\,$}}1
}

\lstdefinestyle{long}{
  numbers=left, numberstyle=\tiny, stepnumber=5, firstnumber=0
}

\newcommand{\codecoq}[1]{\lstinline[language=Coq]{#1}}
\newcommand{\codec}[1]{\lstinline[language=C]{#1}}
\newcommand{\codeacsl}[1]{\lstinline[style=ACSL]{#1}}

\newcommand{\derparun}[2]{\frac{\partial {#1}}{\partial {#2}}}
\newcommand{\derpar}[3]{\frac{\partial^{#1} {#2}}{\partial {#3}^{#1}}}
\newcommand{\eqdef}{{\;\stackrel{\text{def}}{=}\;}}
\newcommand{\floor}[1]{\left\lfloor#1\right\rfloor}

\newcommand{\N}{\mathbb{N}}
\newcommand{\R}{{\mathbb R}}

\newcommand{\Rdeux}{\R^2}
\newcommand{\demi}{\frac{1}{2}}
\newcommand{\prodscal}[3]{\left<#2,#3\right>_{#1}}
\newcommand{\norme}[2]{\left\|#2\right\|_{#1}}
\newcommand{\ps}[2]{\prodscal{}{#1}{#2}}
\newcommand{\n}[1]{\norme{}{#1}}
\newcommand{\nA}[1]{\norme{A(c)}{#1}}
\newcommand{\psdx}[2]{\prodscal{\Delta x}{#1}{#2}}
\newcommand{\ndx}[1]{\norme{\Delta x}{#1}}
\newcommand{\psAh}[2]{\prodscal{A_h(c)}{#1}{#2}}
\newcommand{\nAh}[1]{\norme{A_h(c)}{#1}}

\newcommand{\norm}[1]{\ensuremath{\|#1\|}}

\newcommand{\xmin}{x_\mathrm{min}}
\newcommand{\xmax}{x_\mathrm{max}}
\newcommand{\tmax}{t_\mathrm{max}}
\newcommand{\imax}{i_\mathrm{max}}
\newcommand{\kmax}{k_\mathrm{max}}
\newcommand{\bfx}{{\bf x}}
\newcommand{\bfdeltax}{{\bf \Delta x}}
\newcommand{\eps}{\varepsilon}

\newcommand{\CN}{\mathrm{CN}}

\newcommand{\soft}[1]{#1}
\newcommand{\coq}{\soft{Coq}}
\newcommand{\ssreflect}{\soft{SSReflect}}
\newcommand{\lang}[1]{\soft{#1}}
\newcommand{\ocaml}{\lang{OCaml}}
\newcommand{\clang}{\lang{C}}

\newcommand{\Fost}{F$\!\oint$st}

\newcommand{\ANRCerPanetFost}{This research was supported by the ANR projects
  CerPAN (ANR-05-BLAN-0281-04) and
  {\Fost} (ANR-08-BLAN-0246-01).}

\def\adots{\mathinner{\mkern2mu\raise 1pt\hbox{.}\mkern 3mu\raise 4pt\hbox{.}\mkern1mu\raise 7pt\hbox{{.}}}}

\newenvironment{myproof}{\begin{proof}}{\end{proof}}

\newcommand{\gloss}[1]{\emph{#1}$^\star$}

\newcommand{\apriori}{\emph{a priori}}

\newcommand{\aposteriori}{\emph{a posteriori}}
\newcommand{\eg}{\emph{e.g.}}
\newcommand{\ie}{\emph{i.e.}}

\newcommand{\mysep}{, }
\newcommand{\myappendixreftext}{Appendix~}
\newcommand{\mymultlinenewline}{\null}
\newcommand{\myeqnarraynewline}{\null}

\usepackage[a4paper]{geometry}
\usepackage{RR}
\usepackage{hyperref}
\usepackage{a4wide}

\RRdate{December 2012}

\RRauthor{
  Sylvie Boldo%
  \thanks[toccata]{Projet Toccata.
    {\tt \{Sylvie.Boldo,Jean-Christophe.Filliatre,Guillaume.Melquiond\}%
      @inria.fr}.}%
  \thanks[lri]{LRI, UMR 8623, Universit\'e Paris-Sud, CNRS, Orsay cedex,
    F-91405.}%
  \and Fran\c{c}ois Cl\'ement%
  \thanks[pomdapi]{Projet Pomdapi.
    {\tt \{Francois.Clement,Pierre.Weis\}@inria.fr}.}%
  \and Jean-Christophe Filli\^atre%
  \thanksref{lri}%
  \thanksref{toccata}%
  \and Micaela Mayero%
  \thanks[lipn]{LIPN, UMR 7030, Universit\'e Paris-Nord, CNRS, Villetaneuse,
    F-93430.\goodbreak
    {\tt Micaela.Mayero@lipn.univ-paris13.fr}.}%
  \and Guillaume Melquiond%
  \thanksref{toccata}%
  \thanksref{lri}%
  \and Pierre Weis%
  \thanksref{pomdapi}
}
\authorhead{S. Boldo, F. Cl\'ement, J.-C. Filli\^atre, M. Mayero, G. Melquiond,
  \& P. Weis}

\RRtitle{Calculs en confiance~: une preuve m\'ecanis\'ee de l'\'equation aux
  d\'eriv\'ees partielles au programme effectif}
\RRetitle{Trusting Computations: a Mechanized Proof from Partial Differential
  Equations to Actual Program}
\titlehead{A Mechanized Proof from Partial Differential Equations to Actual
  Program}

\RRnote{\ANRCerPanetFost}

\RRresume{Les programmes informatiques peuvent \'echouer, que ce soit en raison
de comportements exceptionnels, d'acc\`es hors des bornes d'un tableau, ou
simplement d'erreurs de codage.
On ne peut donc pas leur faire confiance aveugl\'ement.
De ce point de vue, les programmes de calcul scientifique ne font pas
exception, et ils pr\'esentent m\^eme des probl\`emes sp\'ecifiques de
pr\'ecision en raison de l'utilisation massive de calculs en virgule
flottante.
Cependant, leur correction est rarement garantie.
En effet, nous avons d\^u \'etendre des m\'ethodes et des outils existants pour
prouver le comportement correct de programmes pour v\'erifier un programme
d'analyse num\'erique existant.
Ce programme~{\clang} impl\'emente le sch\'ema explicite aux diff\'erences
finies centr\'ees du second ordre pour la r\'esolution de l'\'equation des
ondes mono-dimensionnelle.
En fait, nous sommes all\'es beaucoup plus loin puisque nous avons
v\'erifi\'e m\'ecaniquement la convergence du sch\'ema num\'erique pour
obtenir une preuve formelle compl\`ete couvrant tous les aspects de
l'\'equation aux d\'eriv\'ees partielles aux r\'esultats num\'eriques
effectifs.
\`A notre connaissance, c'est la premi\`ere fois qu'une telle preuve
compl\`ete est obtenue.
}
\RRabstract{Computer programs may go wrong due to exceptional behaviors, out-of-bound array
accesses, or simply coding errors.
Thus, they cannot be blindly trusted.
Scientific computing programs make no exception in that respect, and even bring
specific accuracy issues due to their massive use of floating-point
computations.
Yet, it is uncommon to guarantee their correctness.
Indeed, we had to extend existing methods and tools for proving the correct
behavior of programs to verify an existing numerical analysis program.
This {\clang}~program implements the second-order centered finite difference
explicit scheme for solving the 1D wave equation.
In fact, we have gone much further as we have mechanically verified the
convergence of the numerical scheme in order to get a complete formal proof
covering all aspects from partial differential equations to actual numerical
results.
To the best of our knowledge, this is the first time such a comprehensive proof
is achieved.
}

\RRmotcle{\'equation des ondes acoustiques,
preuve formelle d'un programme num\'erique,
convergence d'un sch\'ema num\'erique,
analyse d'erreur d'arrondi.
}
\RRkeyword{Acoustic wave equation\mysep
Formal proof of numerical program\mysep
Convergence of numerical scheme\mysep
Rounding error analysis
}

\RRprojets{Pomdapi et Toccata}

\URRocq

\begin{document}

\RRNo{8197}
\makeRR


\section{Introduction}
\label{sec:introduction}

Given an appropriate set of mathematical equations (such as ODEs or PDEs)
modeling a physical event, the usual simulation process consists of two
stages.
First, the continuous equations are approximated by a set of discrete
equations, called the numerical scheme, which is then proved to be convergent.
Second, the set of discrete equations is implemented as a computer program,
which is eventually run to perform simulations.

The modeling of critical systems requires correctness of the modeling programs
in the sense that there is no runtime error and the computed value is an
accurate solution to the mathematical equations.
The correctness of the program relies on the correctness of the two stages.
Note that we do not consider here the issue of adequacy of the mathematical
equations with the physical phenomenon of interest.
We take the differential equation as a starting point in the simulation
process.
Usually, the discretization stage is justified by a pen-and-paper proof of the
convergence of the selected scheme, while, following~\cite{or:vvs:10}, the
implementation stage is ratified by both code verification and solution
verification.
Code verification (checking for bugs) uses manufactured solutions;
it is called {\em validation by tests} below.
Solution verification (checking for convergence of the numerical scheme at
runtime) usually uses {\aposteriori} error estimates to control the numerical
errors;
it is out of scope for this paper, nevertheless we briefly address the issue in
the final discussion.
The drawback of pen-and-paper proofs is that human beings are fallible and
errors may be left, for example in long and tedious proofs involving a large
number of subcases.
The drawback of validation by tests is that, except for exhaustive testing
which is impossible here, it does not imply a proof of correctness in all
possible cases.
Therefore, one may overestimate the convergence rate, or miss coding errors,
or underestimate round-off errors due to floating-point computations.
In short, this methodology only hints at the correctness of modeling programs
but does not guarantee it.

The fallibility of pen-and-paper proofs and the limitations of validation by
tests is not a new problem, and has been a fundamental concern for a long time
in the computer science community.
The answer to this question came from mathematical logic with the notion of
{\em logical framework} and {\em formal proof}.
A logical framework provides tools to describe mathematical objects and
results, and state theorems to be proved.
Then, the proof of those theorems gets all its logical steps verified in the
logical framework by a computer running a mechanical proof checker.
This kind of proof forbids logical errors and prevents omissions: it is a
formal proof.
Therefore, a formal proof can be considered as a perfect pen-and-paper proof.

Fortunately, logical frameworks also support the definition of computer
programs and the specification of their properties.
The correctness of a program can then be expressed as a formal proof that no
execution of the program will go wrong and that it has the expected
mathematical properties.
A formal proof of a program can be considered as a comprehensive validation by
an exhaustive set of tests.
Note, however, that we verify the program at the source code level and do not
consider here compilation problems, nor external attacks.

Mechanical proof checkers are mainly used to formalize mathematics and are
routinely used to prove programs in the field of integer arithmetic and
symbolic computation.
We apply the same methodology to numerical programs in order to obtain the same
safety level in the scientific computing field.
The simulation process is revisited as follows.
The discretization stage requires some preliminary work in the logical
framework;
we must implement the necessary mathematical concepts and results to describe
continuous and discrete equations (in particular, the notion of convergent
numerical scheme).
Given this mathematical setting, we can write a faithful formal proof of the
convergence of the discrete solution towards the solution to the continuous
problem.
Then, we can specify the modeling program and the properties of the computed
values, and obtain a formal proof of its correctness.
If we specify that computed values are close enough to the actual solution of
the numerical scheme, then the correctness proof of the program ensures the
correctness of the whole simulation process.

This revised simulation process seems easy enough.
However, the difficulty of the necessary formal proofs must not be
underestimated, notably because scientific computing adds specific difficulties
to specifications and proofs.
The discretization stage uses real numbers and real analysis theory.
The usual theorems and tools of real analysis are still in their infancy in
mechanical proof checkers.
In addition, numerical programs use floating-point arithmetic.
Properties of floating-point arithmetic are complex and highly nonintuitive,
which is yet another challenge for formal proofs of numerical programs.

To summarize, the field of scientific computing has the usual difficulties of
formal proof for mathematics and programs, and the specific difficulties of
real analysis and its relationships to floating-point arithmetic.
This complexity explains why mechanical proof checkers are mostly unknown in
scientific computing.
Recent progress~\cite{May01,BolMel11,DinLauMel11,BLM12} in providing mechanical
proof checkers with formalizations of real analysis and IEEE-754 floating-point
arithmetic makes formal proofs of numerical programs tractable nowadays.

In this article, we conduct the formal proof of a very simple {\clang}~program
implementing the second-order centered finite difference explicit scheme for
solving the one-dimensional acoustic wave equation.
This is a first step towards the formal proof of more complex programs used in
critical situations. This article complements a previous publication about the
same experiment~\cite{BCFMMW12}.
This time however, we do not focus on the advances of some formal proof
techniques, but we rather present an overview of how formal methods can be
useful for scientific computing and what it takes to apply them.

Formal proof systems are relatively recent compared with mathematics or
computer science.
The system considered as the first proof assistant is Automath.
It has been designed by de Bruijn in 1967 and has been very influential for the
evolution of proof systems.
As a matter of comparison, the FORTRAN language was born in 1954.
Almost all modern proof assistants then appeared in the 1980s.
In particular, the first version of {\coq} was created in 1984 by Coquand and
Huet.
The ability to reason about numerical programs came much later, as it requires
some formal knowledge of arithmetic and analysis.
In {\coq}, real numbers were formalized in 1999 and floating-point numbers in
2001.
One can note that some of these developments were born from interactions
between several domains, and so is this work.

The formal proofs are too long to be given here \textit{in extenso}, so the
paper only explains general ideas and difficulties.
The annotated {\clang}~program and the full {\coq} sources of the formal
development are available from~\cite{code_and_proofs}.

The paper is organized as follows.
The notion of formal proof and the main formal tools are presented in
Section~\ref{sec:formal_proof}.
Section~\ref{sec:wave} describes the PDE, the numerical scheme, and their
mathematical properties.
Section~\ref{sec:wave_proof} is devoted to the formal proof of the
convergence of the numerical scheme, Section~\ref{sec:round} to the formal
proof of the boundedness on the round-off error, and
Section~\ref{sec:program_proof} to the formal proof of the {\clang}~program
implementing the numerical scheme.
Finally, Section~\ref{sec:discussion} paints a broader picture of the study.

A glossary of terms from the mathematical logic and computer science fields is
given in \myappendixreftext\ref{sec:glossary}.
The main occurrences of \gloss{such terms} are emphasized in the text by using
italic font and superscript star.

\section{Formal Proof}
\label{sec:formal_proof}

Modern mathematics can be seen as the science of abstract objects, {\eg} real
numbers, differential equations.
In contrast, mathematical logic researches the various languages used to define
such abstract objects and reason about them.
Once these languages are formalized, one can manipulate and reason about
mathematical proofs: What makes a valid proof? How can we develop one? And so
on.
This paves the way to two topics we are interested in: mechanical
\gloss{verification} of proofs, and automated deduction of theorems.
In both cases, the use of computer-based tools will be paramount to the success
of the approach.

\subsection{What is a Formal Proof?}
\label{sec:formal_proof1}

When it comes to abstract objects, believing that some properties are true
requires some methods of judgment.
Unfortunately, some of these methods might be fallible: they might be incorrect
in general, or their execution might be lacking in a particular setting.
Logical reasoning aims at eliminating any unjustified assumption and ensuring
that only infallible inferences are used, thus leading to properties that are
believed to be true with the greatest confidence.

The reasoning steps that are applied to deduce from a property believed to be
true a new property believed to be true is called an \gloss{inference rule}.
They are usually handled at a syntactic level: only the form of the statements
matters, their content does not.
For instance, the \emph{modus ponens} rule states that, if both properties
``$A$'' and ``if $A$ then $B$'' hold, then property ``$B$'' holds too, whatever
the meaning of $A$ and $B$.
Conversely, if one deduces ``$B$'' from ``$A$'' and ``if $C$ then $B$'', then
something is amiss: while the result might hold, its proof is definitely
incorrect.

This is where \gloss{formal proofs} show up.
Indeed, since inference rules are simple manipulations of symbols, applying
them or checking that they have been properly applied do not require much
intelligence.
(The intelligence lies in choosing which one to apply.)
Therefore, these tasks can be delegated to a computer running a formal system.
The computer will perform them much more quickly and systematically than a
human being could ever do it.
Assuming that such formal systems have been designed with care,%
\footnote{%
  The core of a formal system is usually a very small program, much smaller
  than any proof it will later have to manipulate, and thus easy to check and
  trust.
  For instance, while expressive enough to tackle any proof of modern
  mathematics, the kernel of HOL Light is just 200 lines long.
}
the results they produce are true with the greatest confidence.

The downside of formal proofs is that they are really low-level;
they are down to the most elementary steps of a reasoning.
It is no longer possible to dismiss some steps of the proofs, trusting the
reader to be intelligent enough to fill the blanks.
Fortunately, since inference rules are mechanical by nature, a formal system
can also try to apply them automatically without any user interaction.
Thus it will produce new results, or at least proofs of known results.
At worst, one could imagine that a formal system applies inference rules
blindly in sequence until a complete proof of a given result is found.
In practice, clever algorithms have been designed to find the proper inference
steps for domain-specific properties.
This considerably eases the process of writing formal proofs.
Note that numerical analysis is not amenable to automatic proving yet, which
means that related properties will require a lot of human interaction, as shown
in Section~\ref{sec:auto}.

It should have become apparent by now that formal systems are primarily aimed
at proving and checking mathematical theorems.
Fortunately, programs can be turned into \gloss{semantically} equivalent
abstract objects that formal systems can manipulate, thus allowing to prove
theorems about programs.
These theorems might be about basic properties of a program, {\eg} it will not
evaluate arrays outside their bounds.
They might also be about higher-level properties, {\eg} the computed results
have such and such properties.
For instance, in this paper, we are interested in proving that the values
computed by the program are actually close to the exact solution to the partial
differential equation.
Note that these verifications are said to be static: they are done once and for
all, yet they cover all the future executions of a program.

Formal verification of a program comes with a disclaimer though, since a
program is not just an abstract object, it also has a concrete behavior once
executed.
Even if one has formally proved that a program always returns the expected
value, mishaps might still happen.
Perfect certainty is unachievable.
First and foremost, the \gloss{specification} of what the program is expected
to compute might be wrong or just incomplete.
For instance, a random generator could be defined as being a function that
takes no input and returns a value between 0 and 1.
One could then formally verify that a given function satisfies such a
specification.
Yet that does not tell anything about the actual randomness of the computed
value: the function might always return the same number while still satisfying
the specification.
This means that formal proofs do not completely remove the need for testing, as
one still needs to make sure specifications are meaningful;
but they considerably reduce the need for exhaustive testing.

Another consideration regarding the extent of confidence in formally verified
programs stems from the fact that programs do not run in isolation, so formal
methods have to make some assumptions.
Basically, they assume that the program executed in the end is the one that was
actually verified and not some variation of it.
This seems an obvious assumption, but practice has shown that a program might
be miscompiled, that some malware might be poking memory at random, that a
computer processor might have design flaws, or even that the electromagnetic
environment might cause bit flips either when verifying the program, or when
executing it.
So the trust in what a program actually computes will still be conditioned to
the trust in the environment it is executed in.
Note that this issue is not specific to verified programs, so they still have
the upper hand over unverified programs.
Moreover, formal methods are also applied to improve the overall trust in a
system: formal verification of hardware design is now routine, and formal
verification of compilers~\cite{Ler09,BJLM13} and operating
systems~\cite{seL410} are bleeding edge research topics.

\subsection{Formal Proof Tools at Work}

There is not a single tool that would allow us to tackle the formal
\gloss{verification} of the {\clang}~program we are interested in.
We will use different tools depending on the kind of abstract objects we want
to manipulate or prove properties about.

The first step lies in running the tool Frama-C over the program
(Section~\ref{sec:tool_fp}).
We have slightly modified the {\clang}~program  by adding comments stating what
the program is expected to compute.
These \gloss{annotations} are just mathematical properties of the program
variables, {\eg} the result variables are close approximations to the values of
the exact solution.
Except for these comments, the code was not modified.
Frama-C takes the program and the annotations and it generates a set of
theorems.
What the tool guarantees is that, if we are able to prove all these theorems,
then the program is formally verified.
Some of these theorems ensure that the execution will not cause exceptional
behaviors: no accesses out of the bounds of the arrays, no overflow during
computations, and so on.
The other theorems ensure that the program satisfies all its annotations.

At this point, we can run tools over the generated theorems, in the hope that
they will automatically find proofs of them.
For instance, Gappa (Section~\ref{sec:tool_gappa}) is suited for proving
theorems stating that floating-point operations do not overflow or that their
round-off error is bounded, while \gloss{SMT solvers}
(Section~\ref{sec:tool_fp}) will tackle theorems stating that arrays are never
accessed out of their bounds.
Unfortunately, more complicated theorems require some user interaction, so we
have used the {\coq} proof assistant (Section~\ref{sec:tool_coq}) to help us in
writing their formal proofs.
This is especially true for theorems that deal with the more
mathematically-oriented aspect of verification, such as convergence of the
numerical scheme.

\subsubsection{{\coq}}
\label{sec:tool_coq}

{\coq}\footnote{\url{http://coq.inria.fr/}} is a formal system that provides
an expressive language to write mathematical definitions, executable
algorithms, and theorems, together with an interactive environment for
proving them~\cite{Coq}.
{\coq}'s formal language combines both a \gloss{higher-order logic} and a
richly-typed \gloss{functional programming} language~\cite{CIC}.
In addition to functions and predicates, {\coq} allows the specification of
mathematical theorems and software \gloss{specifications}, and to interactively
develop formal proofs of those.

The architecture of {\coq} can be split into three parts.
First, there is a relatively small \emph{kernel} that is responsible for
mechanically checking formal proofs.
Given a theorem proved in {\coq}, one does not need to read and understand the
proof to be sure that the theorem statement is correct, one just has to trust
this kernel.

Second, {\coq} provides a proof development system so that the user does not
have to write the low-level proofs that the kernel expects.
There are some interactive proof methods (proof by induction, proof by
contradiction, intermediate lemmas, and so on), some \gloss{decision} and
semi-decision algorithms ({\eg} proving the equality between polynomials), and
a \gloss{tactic} language for letting the user define his or her own proof
methods.
Note that all these high-level proof tools do not have to be trusted, since the
kernel will check the low-level proofs they produce to ensure that all the
theorems are properly proved.

Third, {\coq} comes with a standard library.
It contains a collection of basic and well-known theorems that have already
been formally proved beforehand.
It provides developments and axiomatizations about sets, lists, sorting,
arithmetic, real numbers, and so on.
In this work, we mainly use the Reals standard library~\cite{May01}, which is a
classical axiomatization of an Archimedean ordered complete field.
It comes from the Coq standard library and provides all the basic theorems
about analysis, {\eg} differentials, integrals.
It does not contain more advanced topics such as the Fourier transform and its
properties though.

Here is a short example taken from our \texttt{alpha.v}
file~\cite{code_and_proofs}:
\begin{lstlisting}[language=Coq]
Lemma Rabs_le_trans: forall a b c d : R,
  Rabs (a - c) + Rabs (c - b) <= d -> Rabs (a - b) <= d.
Proof.
intros a b c d H.
replace (a - b) with ((a - c) + (c - b)) by ring.
apply Rle_trans with (2 := H); apply Rabs_triang.
Qed.
\end{lstlisting}
The function \codecoq{Rabs} is the absolute value on real numbers.
The lemma states that, for all real numbers~$a$, $b$, $c$, and~$d$, if
$|a-c|+|c-b| \le d$, then $|a-b| \le d$.
The proof is therefore quite simple.
We first introduce variables and call~\codecoq{H} the hypothesis of the
conditional stating that $|a-c|+|c-b| \le d$.
To prove that $|a-b| \le d$, we first replace $a-b$ with $(a-c)+(c-b)$,
the proof of that being automatic as it is only an algebraic ring equality.
Then, we are left to prove that $|(a-c)+(c-b)| \le d$.
We use transitivity of $\le$, called \codecoq{Rle_trans}, and
hypothesis~\codecoq{H}.
Then, we are left to prove that $|(a-c)+(c-b)| \le |a-c|+|c-b|$.
This is exactly the triangle inequality, called \codecoq{Rabs_triang}.
The proof ends with the keyword \codecoq{Qed}.

The standard library does not come with a formalization of floating-point
numbers.
For that purpose, we use a large {\coq} library called PFF%
\footnote{\url{http://lipforge.ens-lyon.fr/www/pff/}}
initially developed in~\cite{DauRidThe01} and extended with various results
afterwards \cite{Bol04b}.
It is a high-level formalization of the IEEE-754 international standard for
floating-point arithmetic~\cite{ieee-754,Gold91}.
This formalization is convenient for human interactive proofs as testified by
the numerous proofs using it.
The huge number of lemmas available in the library (about 1400) makes it
suitable for a large range of applications.
The library has been superseded since then by the Flocq library~\cite{BolMel11}
and both libraries were used to prove the floating-point results of this work.

\subsubsection{Frama-C, Jessie, Why, and the SMT Solvers}
\label{sec:tool_fp}

We use the Frama-C platform\footnote{\url{http://www.frama-c.cea.fr/}} to
perform formal verification of {\clang}~programs at the source-code level.
Frama-C is an extensible framework that combines \gloss{static analyzers} for
{\clang}~programs, written as plug-ins, within a single tool.
In this work, we use the Jessie plug-in~\cite{marche07plpv} for
\gloss{deductive verification}.
{\clang}~programs are \gloss{annotated} with behavioral contracts written using
the \emph{ANSI {\clang}~Specification Language} (ACSL for short)~\cite{ACSL}.
The Jessie plug-in translates them to the Jessie language~\cite{marche07plpv},
which is part of the Why verification platform~\cite{filliatre07cav}.
This part of the process is responsible for translating the \gloss{semantics}
of~{\clang} into a set of Why logical definitions (to model {\clang}~types,
memory heap, and so on) and Why programs (to model {\clang}~programs).
Finally, the Why platform computes \gloss{verification conditions} from these
programs, using traditional techniques of weakest preconditions~\cite{Dij75},
and emits them to a wide set of existing theorem provers, ranging from
\gloss{interactive proof assistants} to \gloss{automated theorem provers}.
In this work, we use the {\coq} proof assistant (Section~\ref{sec:tool_coq}),
\gloss{SMT solvers} Alt-Ergo~\cite{conchon08entcs}, CVC3~\cite{CVC3}, and
Z3~\cite{Z3}, and the automated theorem prover Gappa
(Section~\ref{sec:tool_gappa}).
Details about automated and interactive proofs can be found in
Section~\ref{sec:auto}.
The dataflow from {\clang}~source code to theorem provers can be depicted as
follows:
\begin{center}
  \includegraphics[width=0.8\textwidth]{verif_c}
\end{center}
More precisely, to run the tools on a {\clang}~program, we use a graphical
interface called gWhy.
A screenshot is displayed in Figure~\ref{fig:screenshot}, in
Section~\ref{sec:program_proof}.
In this interface, we may call one prover on several \gloss{goals}.
We then get a graphical view of how many goals are proved and by which prover.

In ACSL, annotations are written using \gloss{first-order logic}.
Following the \emph{programming by contract} approach, the specifications
involve preconditions, postconditions, and \gloss{loop invariants}.
The contract of the following function states that it computes the square of an
integer \codeacsl{x}, or rather a lower bound on it:
\begin{lstlisting}[style=ACSL]
  //@ requires x >= 0;
  //@ ensures \result * \result <= x;
  int square_root(int x);
\end{lstlisting}
The precondition, introduced with \codeacsl{requires}, states that the
argument~\codeacsl{x} is nonnegative.
Whenever this function is called, the toolchain will generate a theorem stating
that the input is nonnegative.
The user then has to prove it to ensure the program is correct.
The postcondition, introduced with \codeacsl{ensures}, states the property
satisfied by the return value \lstinline[style=ACSL]{\result}.
An important point is that, in the specification, arithmetic operations are
mathematical, not machine operations.
In particular, the product \lstinline[style=ACSL]{\result * \result} cannot
overflow.
Simply speaking, we can say that {\clang}~integers are reflected within
specifications as mathematical integers, in an obvious way.

The translation of floating-point numbers is more subtle, since one needs to
talk about both the value actually computed by an expression, and the ideal
value that would have been computed if we had computers able to work on real
numbers.
For instance, the following excerpt from our {\clang}~program specifies the
relative error on the content of the \codeacsl{dx} variable, which represents
the grid step $\Delta x$ (see Section~\ref{sec:discrete}):
\begin{lstlisting}[style=ACSL]
  dx = 1./ni;
  /*@ assert
    @   dx > 0. && dx <= 0.5 &&
    @   \abs(\exact(dx) - dx) / dx <= 0x1.p-53;
    @ */
\end{lstlisting}
The identifier \codeacsl{dx} represents the value actually computed (seen as a
real number), while the expression \lstinline[style=ACSL]{\exact(dx)}
represents the value that would have been computed if mathematical operators
had been used instead of floating-point operators.
Note that \codeacsl{0x1.p-53} is a valid ACSL literal (and {\clang}~too)
meaning $1 \cdot 2^{-53}$ (which is also the machine epsilon on {\sf binary64}
numbers).

\subsubsection{Gappa}
\label{sec:tool_gappa}

The Gappa tool\footnote{\url{http://gappa.gforge.inria.fr/}} adapts the
\gloss{interval-arithmetic} paradigm to the proof of properties that occur when
verifying numerical applications~\cite{DauMel10}.
The inputs are logical formulas quantified over real numbers whose
\gloss{atoms} are usually enclosures of arithmetic expressions inside numeric
intervals.
Gappa answers whether it succeeded in verifying it.
In order to support program \gloss{verification}, one can use \emph{rounding}
functions inside expressions.
These unary operators take a real number and return the closest real number in
a given direction that is representable in a given binary floating-point
format.
For instance, assuming that operator~$\circ$ rounds to the nearest
{\sf binary64} floating-point number, the following formula states that the
relative error of the floating-point addition is bounded~\cite{Gold91}:
\[\forall x,y\in\R,~\exists\eps\in\R,~|\eps| \le 2^{-53} \mbox{ and }
{\circ}(\circ(x)+\circ(y)) = (\circ(x) + \circ(y)) (1 + \eps).\]

Converting \gloss{straight-line} numerical programs to Gappa logical formulas
is easy and the user can provide additional hints if the tool were to fail to
verify a property.
The tool is specially designed to handle codes that perform convoluted
floating-point manipulations.
For instance, it has been successfully used to verify a state-of-the-art
library of correctly-rounded elementary functions~\cite{DinLauMel11}.
In the current work, Gappa has been used to check much simpler properties.
In particular, no user hint was needed to automatically prove them.
Yet the length of their proofs would discourage even the most dedicated users
if they were to be manually handled.
One of the properties is the round-off error of a local evaluation of the
numerical scheme (Section~\ref{sec:local-round-off-error}).
Other properties mainly deal with proving that no exceptional behavior occurs
while executing the program: due to the initial values, all the computed values
are sufficiently small to never cause overflow.

Verification of some formulas requires reasonings that are so long and
intricate~\cite{DinLauMel11}, that it might cast some doubts on whether an
automatic tool can actually succeed in proving them.
This is especially true when the tool ends up proving a property stronger than
what the user expected.
That is why Gappa also generates a formal proof that can be mechanically
checked by a proof assistant.
This feature has served as the basis for a {\coq} \gloss{tactic} for
automatically proving \gloss{goals} related to floating-point and real
arithmetic~\cite{BolFilMel09}.
Note that Gappa itself is not verified, but since {\coq} verifies the proofs
that Gappa generates, the goals are formally proved.

This tactic has been used whenever a \gloss{verification condition} would have
been directly proved by Gappa, if not for some confusing notations or encodings
of matrix elements.
We just had to apply a few basic {\coq} tactics to put the goal into the proper
form and then call the Gappa tactic to prove it automatically.

\section{Numerical Scheme for the Wave Equation}
\label{sec:wave}

We have chosen to study the numerical solution to the one-dimensional
acoustic wave equation using the second-order centered explicit scheme as
it is simple, yet representative of a wide class of scientific computing
problems.
First, following~\cite{bec:esn:09}, we describe and state the different notions
necessary for the implementation of the numerical scheme and its analysis.
Then, we present the annotations added in the source code to specify the
behavior of the program.

\subsection{Continuous Equation}
\label{sec:continuous}

We consider $\Omega=[\xmin,\xmax]$, a one-dimensional homogeneous acoustic
medium characterized by the constant propagation velocity~$c$.
Let~$p(x,t)$ be the acoustic quantity, {\eg} the transverse displacement
of a vibrating string, or the acoustic pressure.
Let~$p_0(x)$ and~$p_1(x)$ be the initial conditions.
Let us consider homogeneous Dirichlet boundary conditions.

The one-dimensional acoustic problem on~$\Omega$ is set by
\begin{eqnarray}
  \label{eq:L}
  \forall t \ge 0,\; \forall x \in \Omega, & &
  (L (c) \, p) (x, t) \eqdef
  \derpar{2}{p}{t} (x, t) + A (c) \, p (x, t) = 0, \\
  \label{eq:L1}
  \forall x \in \Omega, & &
  (L_1 \, p) (x, 0) \eqdef
  \derparun{p}{t} (x, 0) =
  p_1 (x), \\
  \label{eq:L0}
  \forall x \in \Omega, & &
  (L_0 \, p) (x, 0) \eqdef
  p (x, 0) =
  p_0 (x), \\
  \label{eq:dir}
  \forall t \ge 0, & &
  p (\xmin, t) =
  p (\xmax, t) =
  0
\end{eqnarray}
where the differential operator~$A(c)$ acting on function~$q$ is defined as
\begin{equation}
  \label{eq:A}
  A (c) q \eqdef - c^2 \derpar{2}{q}{x}.
\end{equation}

We assume that under reasonable regularity conditions on the Cauchy data~$p_0$
and~$p_1$, for each $c>0$, there exists a unique solution~$p$ to the
initial-boundary value problem defined by Equations~(\ref{eq:L})
to~(\ref{eq:A}).
Of course, it is well-known that the solution to this partial differential
equation is given by d'Alembert's formula~\cite{dal:rcf:47}.
But simply assuming existence of a solution instead of exhibiting it ensures
that our approach scales to more general cases for which there is no known
analytic expression of a solution, {\eg} when propagation velocity~$c$
depends on space variable~$x$.

We introduce the positive definite quadratic quantity
\begin{equation}
  \label{eq:energiecontinue}
  E (c) (p) (t) \eqdef
  \demi \n{\derparun{p}{t} (\cdot, t)}^2 +
  \demi \nA{p (\cdot, t)}^2
\end{equation}
where $\ps{q}{r}\eqdef\int_\Omega q(x)r(x)dx$,
$\n{q}^2\eqdef\ps{q}{q}$ and $\nA{q}^2\eqdef\ps{A(c)\,q}{q}$.
The first term is interpreted as the kinetic energy, and the second term as the
potential energy, making~$E$ the mechanical energy of the acoustic system.

Let~$\tilde{p}_0$ (resp.~$\tilde{p}_1$) represent the function defined on the
entire real axis~$\R$ obtained by successive antisymmetric extensions in space
of~$p_0$ (resp.~$p_1$).
For example, we have, for all $x\in[2\xmin-\xmax,\xmin]$,
$\tilde{p_0}(x)=-p_0(2\xmin-x)$.
The image theory~\cite{joh:pde:86} stipulates that the solution of the wave
equation defined by Equations~(\ref{eq:L}) to~(\ref{eq:A}) coincides on
domain~$\Omega$ with the solution of the same wave equation but set on the
entire real axis~$\R$, without the Dirichlet boundary
condition~(\ref{eq:dir}), and with extended Cauchy data~$\tilde{p}_0$
and~$\tilde{p}_1$.

\subsection{Discrete Equations}
\label{sec:discrete}

Let us consider the time interval $[0,\tmax]$.
Let~$\imax$ (resp.~$\kmax$) be the number of intervals of the space
(resp. time) discretization.
We define%
\footnote{{\em Floor} notation $\floor{.}$ denotes rounding to
  an integer towards minus infinity.}
\begin{eqnarray}
  \label{eq:dx}
  \Delta x \eqdef \frac{\xmax - \xmin}{\imax}, & \quad &
  i_{\Delta x} (x) \eqdef \floor{\frac{x - \xmin}{\Delta x}}, \\
  \label{eq:dt}
  \Delta t \eqdef \frac{\tmax}{\kmax}, & \quad &
  k_{\Delta t} (t) \eqdef \floor{\frac{t}{\Delta t}}.
\end{eqnarray}
The regular discrete grid approximating $\Omega\times[0,\tmax]$ is
defined by\footnote{For integers~$n$ and~$p$, the notation $[n..p]$
  denotes the integer range $[n,p]\cap\N$.}
\begin{equation}
  \label{eq:xik}
  \forall k \in [0..\kmax],\,
  \forall i \in [0..\imax], \quad
  \bfx_i^k \eqdef
  (x_i, t^k) \eqdef
  (\xmin + i \Delta x, k \Delta t).
\end{equation}

For a function~$q$ defined over $\Omega\times[0,\tmax]$ (resp. $\Omega$), the
notation~$q_{\rm h}$ (with a roman index ${\rm h}$) denotes any discrete
approximation of~$q$ at the points of the grid, {\ie} a discrete function
over $[0..\imax]\times[0..\kmax]$ (resp. $[0..\imax]$).
By extension, the notation~$q_{\rm h}$ is also a shortcut to denote the matrix
$(q_i^k)_{0\leq i\leq \imax,0\leq k\leq\kmax}$
(resp. the vector $(q_i)_{0\leq i\leq \imax}$).
The notation~$\bar{q}_{\rm h}$ (with a bar over it) is reserved to
specify the values of~$q_{\rm h}$ at the grid points,
$\bar{q}_i^k\eqdef q(\bfx_i^k)$ (resp. $\bar{q}_i\eqdef q(x_i)$).

Let~$u_{0,{\rm h}}$ and~$u_{1,{\rm h}}$ be two discrete functions over
$[0..\imax]$.
Let~$s_{\rm h}$ be a discrete function over $[0..\imax]\times[0..\kmax]$.
Then, the discrete function~$p_{\rm h}$ over $[0..\imax]\times[0..\kmax]$ is
the
{\em solution of the second-order centered finite difference explicit scheme},
when the following set of equations holds:
\begin{multline}
  \label{eq:Lh}
  \forall k \in [2..\kmax],\,
  \forall i \in [1..\imax-1], \\
  (L_{\rm h} (c) \, p_{\rm h})_i^k \eqdef
  \frac{p_i^k - 2 p_i^{k - 1} + p_i^{k - 2}}{\Delta t^2} +
  (A_{\rm h} (c) \, p_{\rm h}^{k - 1})_i =
  s_i^{k - 1},
\end{multline}
\begin{eqnarray}
  \label{eq:L1h}
  \forall i \in [1..\imax-1], & &
  (L_{1,{\rm h}} (c) \, p_{\rm h})_i \eqdef
  \frac{p_i^1 - p_i^0}{\Delta t} +
  \frac{\Delta t}{2} (A_{\rm h} (c) \, p_{\rm h}^0)_i =
  u_{1,i}, \\
  \label{eq:L0h}
  \forall i \in [1..\imax-1], & &
  (L_{0,{\rm h}} \, p_{\rm h})_i \eqdef
  p_i^0 =
  u_{0,i}, \\
  \label{eq:dirh}
  \forall k \in [0..\kmax], & &
  p_0^k =
  p_{\imax}^k =
  0
\end{eqnarray}
where the matrix~$A_{\rm h}(c)$ (discrete analog of~$A(c)$) is defined on
vectors~$q_{\rm h}$ by
\begin{equation}
  \label{eq:Ah}
  \forall i \in [1..\imax-1], \quad
  \left( A_{\rm h} (c) \, q_{\rm h} \right)_i \eqdef
  -c^2 \frac{q_{i+1} - 2 q_i + q_{i-1}}{\Delta x^2}.
\end{equation}
Note the use of a second-order approximation of the first derivative in time
in Equation~(\ref{eq:L1h}).

A discrete analog of the energy is also defined by
\begin{equation}
  \label{eq:discreteenergy}
  E_{\rm h} (c) (p_{\rm h}) ^ {k+\demi} \eqdef
  \demi \ndx{\frac{p_{\rm h}^{k+1} - p_{\rm h}^k}{\Delta t}}^2 +
  \demi \psAh{p_{\rm h}^k}{p_{\rm h}^{k+1}}
\end{equation}
where, for any vectors~$q_{\rm h}$ and~$r_{\rm h}$,
\[
\begin{array}{ll}
  \psdx{q_{\rm h}}{r_{\rm h}} \eqdef
  \displaystyle \sum_{i=1}^{\imax-1} q_ir_i\Delta x, &
  \ndx{q_{\rm h}}^2 \eqdef \psdx{q_{\rm h}}{q_{\rm h}}, \\
  \psAh{q_{\rm h}}{r_{\rm h}} \eqdef
  \psdx{A_{\rm h}(c)\,q_{\rm h}}{r_{\rm h}}, &
  \nAh{q_{\rm h}}^2 \eqdef \psAh{q_{\rm h}}{q_{\rm h}}.
\end{array}
\]
Note that $\nAh{.}$ is a semi-norm.
Thus, the discrete energy~$E_h(c)$ is targeted to be a positive semidefinite
quadratic form.

Note that the numerical scheme is parameterized by the discrete Cauchy
data~$u_{0,{\rm h}}$ and~$u_{1,{\rm h}}$, and by the discrete source
term~$s_{\rm h}$.
When these input data are respectively approximations of the exact
functions~$p_0$ and~$p_1$ ({\eg} when
$u_{0,{\rm h}}=\bar{p_0}_{\rm h}$, $u_{1,{\rm h}}=\bar{p_1}_{\rm h}$, and
$s_{\rm h}\equiv 0$), the discrete solution~$p_{\rm h}$ is an approximation of
the exact solution~$p$.

The remark about image theory holds here too: we may replace the use of
Dirichlet boundary conditions in Equation~(\ref{eq:dirh}) by considering
extended discrete Cauchy data~${\bar{\tilde{p_0}}}_{\rm h}$
and~${\bar{\tilde{p_1}}}_{\rm h}$.

Note also that, as well as in the continuous case when a source term is
considered on the right-hand side of Equation~(\ref{eq:L}), the discrete
solution of the numerical scheme defined by Equations~(\ref{eq:Lh})
to~(\ref{eq:Ah}) can be obtained by the discrete space--time convolution of the
discrete fundamental solution and the discrete source term.
This will be useful in Sections~\ref{sec:global-round-off-error}
and~\ref{sec:future_work}.

\subsection{Convergence}
\label{sec:wave_conv}

The main properties required for a numerical scheme are the convergence, the
consistency, and the stability.
A numerical scheme is {\em convergent} when the convergence error,
{\ie} the difference between exact and approximated solutions, tends to zero
with respect to the discretization parameters.
It is {\em consistent with the continuous equations} when the truncation
error, {\ie} the difference between exact and approximated equations, tends to
zero with respect to the discretization parameters.
It is {\em stable} if the approximated solution is bounded when
discretization parameters tend to zero.

The Lax equivalence theorem stipulates that consistency implies equivalence
between stability and convergence, {\eg} see~\cite{str:fds:89}.
Consistency proof is usually straightforward, and stability proof is typically
much easier than convergence proof.
Therefore, in practice, the convergence of numerical schemes is obtained by
using the Lax equivalence theorem once consistency and stability are
established.
Unfortunately, we cannot follow this popular path, since building a
\gloss{formal proof} of the Lax equivalence theorem is quite challenging.
Assuming such a powerful theorem is not an option either, since it would
jeopardize the whole formal proof.
Instead, we formally prove that consistency and stability imply convergence in
the particular case of the second-order centered scheme for the wave equation.

The Fourier transform is a popular tool to study the convergence of numerical
schemes.
Unfortunately, the formalization of the Lebesgue integral theory and the
Fourier transform theory does not yet exist in {\coq}
(Section~\ref{sec:tool_coq});
such a development should not encounter any major difficulty, except for its
human cost.
As an alternative, we consider energy-based techniques.
The stability is then expressed in terms of a discrete energy which only
involves finite summations (to compute discrete dot products).
Energy-based approaches are less precise because they follow {\apriori} error
estimates; but, unlike the Fourier analysis approach, they can be extended to
the heterogeneous case, and to non uniform grids.

The CFL condition (for Courant-Friedrichs-Lewy, see~\cite{cfl:pde:67}) states
for the one-dimensional acoustic wave equation that the Courant number
\begin{equation}
  \label{eq:cn}
  \CN \eqdef \frac{c \Delta t}{\Delta x}
\end{equation}
should not be greater than~1.
To simplify the formal proofs, we use a more restrictive CFL condition.
First, we rule out the particular case $\CN=1$ since it is not so important in
practice, and would imply a devoted structure of proofs.
Second, the convergence proof based on energy techniques requires~$\CN$ to stay
away from~1 (for instance, constant~$C_2$ in
Equation~(\ref{eq:energystabilityconstants}) may explode when $\xi$ tends
to~0).
Thus, we parameterize the CFL condition with a real number~$\xi\in]0,1[$:
\begin{equation}
  \label{eq:cfl}
  \mathrm{CFL} (\xi) \eqdef \CN \leq 1 - \xi.
\end{equation}
Our formal proofs are valid for all~$\xi$ in $]0,1[$.
However, to deal with floating-point arithmetic, the {\clang}~code is annotated
with values of~$\xi$ down to $2^{-50}$.

For the numerical scheme defined by Equations~(\ref{eq:Lh}) to~(\ref{eq:Ah}),
the convergence error~$e_{\rm h}$ and the truncation error~$\eps_{\rm h}$ are
defined by
\begin{eqnarray}
  \label{eq:conv_error}
  \forall k \in [0..\kmax],\,
  \forall i \in [0..\imax], \quad
  e_i^k & \eqdef & \bar{p}_i^k - p_i^k, \\
  \label{eq:trunc_error_k}
  \forall k \in [2..\kmax],\,
  \forall i \in [1..\imax-1], \quad
  \eps_i^k & \eqdef & (L_{\rm h} (c) \, \bar{p}_{\rm h})_i^k, \\
  \label{eq:trunc_error_1}
  \forall i \in [1..\imax-1], \quad
  \eps_i^1 & \eqdef & (L_{1,{\rm h}} (c) \, \bar{p}_{\rm h})_i - \bar{p}_{1,i}, \\
  \label{eq:trunc_error_0}
  \forall i \in [1..\imax-1], \quad
  \eps_i^0 & \eqdef & (L_{0,{\rm h}} \bar{p}_{\rm h})_i - \bar{p}_{0,i}, \\
  \forall k \in [0..\kmax], \quad
  \eps_0^k = \eps_{\imax}^k & \eqdef & 0.
\end{eqnarray}
Note that, when the input data of the numerical scheme for the approximation of
the exact solution are given by $u_{0,{\rm h}}=\bar{p_0}_{\rm h}$,
$u_{1,{\rm h}}=\bar{p_1}_{\rm h}$, and $s_{\rm h}\equiv 0$,
then the convergence error~$e_{\rm h}$ is itself solution of the same numerical
scheme with discrete inputs corresponding to the truncation error:
$u_{0,{\rm h}}=\eps_{\rm h}^0=0$, $u_{1,{\rm h}}=\eps_{\rm h}^1$, and
$s_{\rm h}=\left(k\mapsto\eps_{\rm h}^{k+1}\right)$.
This will be useful in Section~\ref{sec:conv}.

In Section~\ref{sec:conv}, we discuss the formal proof that the numerical
scheme is {\em convergent of order~($m$, $n$) uniformly on the interval
  $[0,\tmax]$} if the convergence error satisfies%
\footnote{%
  The big~$O$ notation is defined in Section~\ref{sec:o}.
  The function $k_{\Delta t}$ is defined in Equation~(\ref{eq:dt}).}
\begin{equation}
  \label{eq:convergence}
  \ndx{e_{\rm h}^{k_{\Delta t}(t)}} = O_{[0,\tmax]} (\Delta x^m + \Delta t^n).
\end{equation}

In Section~\ref{sec:consist}, we discuss the formal proof that the numerical
scheme is {\em consistent with the continuous problem at order~($m$, $n$)
  uniformly on interval $[0,\tmax]$} if the truncation error satisfies
\begin{equation}
  \label{eq:consistency}
  \ndx{\eps_{\rm h}^{k_{\Delta t}(t)}} = O_{[0,\tmax]} (\Delta x^m + \Delta t^n).
\end{equation}

In Section~\ref{sec:stab}, we discuss the formal proof that the numerical
scheme is {\em energetically stable uniformly on interval $[0,\tmax]$} if the
discrete energy defined by Equation~(\ref{eq:discreteenergy}) satisfies
\begin{multline}
  \label{eq:energystability}
  \exists \alpha, C_1, C_2 > 0,\,
  \forall t \in [0, \tmax],\,
  \forall \Delta x, \Delta t > 0, \quad
  \sqrt{\Delta x^2 + \Delta t^2} < \alpha \Rightarrow \\
  \sqrt{E_{\rm h} (c) (p_{\rm h})^{k_{\Delta t}(t) + \demi}} \le
  C_1 + C_2 \Delta t \sum_{k^\prime=1}^{k_{\Delta t}(t)}
  \norme{\Delta x}{\left(i \mapsto s_i^{k^\prime}\right)}.
\end{multline}
Note that constants~$C_1$ and~$C_2$ may depend on the discrete Cauchy
data~$u_{0,{\rm h}}$ and~$u_{1,{\rm h}}$.

\subsection{{\clang}~Program}
\label{sec:program}

\begin{lstlisting}[
    float, label={l:unannotated_code},
    caption={The main part of the {\clang}~code, without annotations.},
    language=C, style=long
  ]
/* Compute the constant coefficient of the stiffness matrix. */
a1 = dt/dx*v;
a  = a1*a1;

/* First initial condition and boundary conditions. */
/* Left boundary. */
p[0][0] = 0.;
/* Time iteration -1 = space loop. */
for (i=1; i<ni; i++) {
  p[i][0] = p0(i*dx);
}
/* Right boundary. */
p[ni][0] = 0.;

/* Second initial condition (with p1=0) and boundary conditions. */
/* Left boundary. */
p[0][1] = 0.;
/* Time iteration 0 = space loop. */
for (i=1; i<ni; i++) {
  dp = p[i+1][0] - 2.*p[i][0] + p[i-1][0];
  p[i][1] = p[i][0] + 0.5*a*dp;
}
/* Right boundary. */
p[ni][1] = 0.;

/* Evolution problem and boundary conditions. */
/* Propagation = time loop. */
for (k=1; k<nk; k++) {
  /* Left boundary. */
  p[0][k+1] = 0.;
  /* Time iteration k = space loop. */
  for (i=1; i<ni; i++) {
    dp = p[i+1][k] - 2.*p[i][k] + p[i-1][k];
    p[i][k+1] = 2.*p[i][k] - p[i][k-1] + a*dp;
  }
  /* Right boundary. */
  p[ni][k+1] = 0.;
}
\end{lstlisting}

We assume that $\xmin=0$, $\xmax=1$, and that the absolute value of the exact
solution is bounded by~1.
The main part of the {\clang}~program is listed in
Listing~\ref{l:unannotated_code}.
Grid steps $\Delta x$ and~$\Delta t$ are respectively represented by
variables~\codec{dx} and~\codec{dt}, grid sizes~$\imax$ and~$\kmax$ by
variables~\codec{ni} and~\codec{nk}, and the propagation velocity~$c$ by the
variable~\codec{v}.
The Courant number~$\CN$ is represented by variable~\codec{a1}.
The initial value~$u_{0,{\rm h}}$ is represented by the function~\codec{p0}.
Other input data~$u_{1,{\rm h}}$ and~$s_{\rm h}$ are supposed to be
zero and are not represented.
The discrete solution~$p_{\rm h}$ is represented by the two-dimensional
array~\codec{p} of size $(\imax+1)(\kmax+1)$.
Note that this is a simple implementation.
A more efficient one would only store two time steps.

\subsection{Program Annotations}
\label{sec:annot}

As explained in Section~\ref{sec:tool_fp}, the {\clang}~code is enriched with
\gloss{annotations} that specify the behavior of the program: what it requires
on its inputs and what it ensures on its outputs.
We describe here the chosen specification for this program.
The full annotations are given in
\myappendixreftext\ref{sec:annotated_source}.

The annotations can be separated into two sets.
The first one corresponds to the mathematics: the exact solution of the wave
equation and its properties.
It defines the required values (the exact solution $p$, and its initialization
$p_0$).
It also defines the derivatives of $p$: \codeacsl{psol_1}, \codeacsl{psol_2},
\codeacsl{psol_11}, and \codeacsl{psol_22} respectively stand for
$\derparun{p}{x}$, $\derparun{p}{t}$, $\derpar{2}{p}{x}$, and
$\derpar{2}{p}{t}$.
The value of the derivative of~$f$ at point~$x$ is defined as the limit of
$\frac{f(x+h)-f(x)}{h}$ when $h \rightarrow 0$.
As the ACSL annotations are only \gloss{first-order}, these definitions are
quite cumbersome: each derivative needs 5 lines to be defined.
Here is the example of \codeacsl{psol_1}, {\ie} $\derparun{p}{x}$:
\begin{lstlisting}[style=ACSL]
/*@ logic real psol_1(real x, real t);
  @ axiom psol_1_def:
  @   \forall real x; \forall real t;
  @   \forall real eps; \exists real C; 0 < C && \forall real dx;
  @   0 < eps ==> \abs(dx) < C ==>
  @   \abs((psol(x + dx, t) - psol(x, t)) / dx - psol_1(x, t)) < eps;
  @ */
\end{lstlisting}
Note the different treatment of the positiveness for the existential
variable~\codeacsl{C}, and for the free variable~\codeacsl{eps}.

We also assume that the solution actually solves Equations~(\ref{eq:L})
to~(\ref{eq:A}).
The last property needed on the exact solution is its regularity.
We require it to be close to its Taylor approximations of degrees 3 and 4 on
the whole space interval (see Section~\ref{sec:wave_proof}).
For instance, the following annotation states the property for degree~3:
\begin{lstlisting}[style=ACSL]
/*@ axiom psol_suff_regular_3:
  @   0 < alpha_3 && 0 < C_3 &&
  @   \forall real x; \forall real t; \forall real dx; \forall real dt;
  @   0 <= x <= 1 ==> \sqrt(dx * dx + dt * dt) <= alpha_3 ==>
  @   \abs(psol(x + dx, t + dt) - psol_Taylor_3(x, t, dx, dt)) <=
  @     C_3 * \abs(\pow(\sqrt(dx * dx + dt * dt), 3));
  @ */
\end{lstlisting}

The second set of annotations corresponds to the properties and
\gloss{loop invariant} needed by the program.
For example, we require the matrix to be \textit{separated}.
This means that a line of the matrix should not overlap with another line.
Otherwise a modification would alter another entry in the matrix.
The predicate \codeacsl{analytic_error} that is used as a loop invariant is
declared in the annotations and defined in the {\coq} files.

The initialization functions are only specified: there is no {\clang}~code
given and no proof.
We assume a reasonable behavior for these functions.
More precisely, this corresponds first to the function \codec{array2d_alloc}
that initializes the matrix and \codec{p_zero} that produces an approximation
of the $p_0$ function.
The use of \gloss{deductive verification} allows our proof to be generic with
respect to $p_0$ and its implementation.

The preconditions of the main function are as follows:
\begin{itemize}
\item $\imax$ and~$\kmax$ must be greater than one, but small enough so that
  $\imax+1$ and  $\kmax+1$ do not overflow;
\item the $\mathrm{CFL}(2^{-50})$ condition must be satisfied (see
  Equations~(\ref{eq:cfl}) and~(\ref{eq:cn}));
\item the grid sizes~$\Delta x$ and~$\Delta t$ must be small enough to ensure
  the convergence of the scheme (the precise value is given in
  Section~\ref{sec:total_error});
\item the floating-point values computed for the grid sizes must be close to
  their mathematical values;
\item to prevent exceptional behavior in the computation of $a$, $\Delta t$
 must be greater than $2^{-1000}$ and $\frac{c \Delta t}{\Delta x}$ must be
 greater than $2^{-500}$.
\end{itemize}
The last hypothesis gets very close to the underflow threshold: the smallest
positive floating-point number is a subnormal of value~$2^{-1074}$.

There are two postconditions, corresponding to a bound on the convergence error
and a bound on the round-off error.
See Sections~\ref{sec:wave_proof} and~\ref{sec:round} for more details.

\section{Formal Proof of Convergence of the Scheme}
\label{sec:wave_proof}

In~\cite{BoldoFilliatre07}, the method error is the distance between the
discrete value computed with exact real arithmetic and the ideal mathematical
value, {\ie} the exact solution;
here, it reduces to the convergence error.
Round-off errors, due to the use of floating-point arithmetic, are handled in
Section~\ref{sec:round}.
First, we present the notions necessary to the
\gloss{formal specification and proof} of boundedness of the method error.
Then, taking inspiration from~\cite{bec:esn:09}, we sketch the main steps of
the proofs of consistency, stability, and convergence of the numerical scheme,
and we point out some tricky aspects that may not be visible in pen-and-paper
proofs.
Full formal proofs are available from~\cite{code_and_proofs}.
Feedback about issues showing up when formalizing a pen-and-paper proof is
discussed in Section~\ref{sec:mathsandproofs}.

\subsection{Big~O, Taylor Approximations, and Regularity}
\label{sec:o}

Proving the consistency of the numerical scheme with the continuous equations
requires some assumptions on the regularity of the exact solution.
This regularity is expressed as the existence of Taylor approximations of the
exact solution up to some appropriate order.
For the sake of simplicity, we chose to prove uniform consistency of the scheme
(see definitions at the end of Section~\ref{sec:wave_conv});
this implies the use of uniform Taylor approximations.

These uniform approximations involve {\em uniform big~$O$} equalities of the
form
\begin{equation}
  \label{eq:bigO}
  F (\bfx, \bfdeltax) = O_{{\Omega_\bfx},\Omega_\bfdeltax} (g(\bfdeltax)).
\end{equation}
Variable~$\bfx$ stands for points around which Taylor approximations are
considered, and variable~$\bfdeltax$ stands for a change in~$\bfx$ in the
neighborhood of zero.
Note that we explicitly state the set~$\Omega_\bfx$ on which the approximation
is uniform.
As for the set~$\Omega_\bfdeltax$, we will later omit it when it spans the full
space, here~$\Rdeux$.
Function~$F$ goes from $\Omega_\bfx\times\Omega_\bfdeltax$ to~$\R$, and~$g$
goes from~$\Omega_\bfdeltax$ to~$\R$.
As all equalities embedding big~O notations, Equation~(\ref{eq:bigO}) actually
means that function $\bfdeltax\mapsto F(\bfx,\bfdeltax)$ belongs to the set
$O_{\Omega_\bfx}(g(\bfdeltax))$ of functions of~$\bfdeltax$ with the same
asymptotic growth rate as~$g$, uniformly for all~$\bfx$ in~$\Omega_\bfx$.

The formal definition behind Equation~(\ref{eq:bigO}) is
\begin{equation}
  \label{eq:Oups}
  \exists \alpha, C > 0,\,
  \forall \bfx \in \Omega_\bfx,\, \forall \bfdeltax \in \Omega_\bfdeltax,
  \norm{\bfdeltax} < \alpha \Rightarrow
  | F (\bfx, \bfdeltax) | \le C |g(\bfdeltax)|.
\end{equation}
Its translation in {\coq} is straightforward.%
\footnote{%
  \codecoq{Rabs} is the absolute value of real numbers, and
  \codecoq{norm_l2} is the $L^2$ norm on~$\Rdeux$.
}
\begin{lstlisting}[language=Coq]
Definition
    OuP (A : Type) (P : R * R -> Prop) (f : A -> R * R -> R) (g : R * R -> R) :=
  exists alpha : R, exists C : R, 0 < alpha /\ 0 < C /\
  forall X : A, forall dX : R * R, P dX -> norm_l2 dX < alpha ->
  Rabs (f X dX) <= C * Rabs (g dX).
\end{lstlisting}
Argument type~\codecoq{A} stands for domain~$\Omega_\bfx$, allowing to deal
with either a subset of~$\Rdeux$ or of~$\R$.
Domain~$\Omega_\bfdeltax$ is formalized by its indicator function on~$\Rdeux$,
the argument predicate~\codecoq{P} (\codecoq{Prop} is the type of logical
propositions).

Let~$\bfx=(x,t)$ be a point in~$\Rdeux$,
let~$f$ be a function from~$\Rdeux$ to~$\R$, and
let~$n$ be a natural number.
The Taylor polynomial of order~$n$ of~$f$ at point~$\bfx$ is the function
$\mathrm{TP}_n(f,\bfx)$ defined over~$\Rdeux$ by
\[
\mathrm{TP}_n (f, \bfx) (\Delta x, \Delta t) \eqdef
\sum_{p=0}^n
\frac{1}{p!} \left(
  \sum_{m=0}^p \binom{p}{m}
  \frac{\partial^p f}{\partial x^m \partial t^{p-m}} (\bfx)
  \cdot \Delta x^m \cdot \Delta t^{p-m}
\right).
\]

Let~$\Omega_\bfx$ be a subset of~$\R^2$.
The Taylor polynomial is a
{\em uniform approximation of order~$n$ of~$f$ on~$\Omega_\bfx$}
when the following uniform big~O equality holds:%
\footnote{Domain~$\Omega_\bfdeltax$ is implicitly considered to be~$\Rdeux$.}
\[
f (\bfx + \bfdeltax) - \mathrm{TP}_n (f, \bfx) (\bfdeltax) =
O_{\Omega_\bfx} \left(\norm{\bfdeltax}^{n+1}\right).
\]
Function~$f$ is then
{\em sufficiently regular of order~$n$ uniformly on~$\Omega_\bfx$}
when all its Taylor polynomials of order smaller than~$n$ are uniform
approximations on~$\Omega_\bfx$.

Of course, we could have expressed regularity requirements by using the usual
notion of differentiability class.
Indeed, it is well-known that functions of class~$C^{n+1}$ over some compact
domain are actually sufficiently regular of order~$n$ uniformly on that
domain.
But, besides the suboptimality of the statement, we did not want to formalize
all the technical details of its proof.

\subsection{Consistency}
\label{sec:consist}

The numerical scheme defined by Equations~(\ref{eq:Lh}) to~(\ref{eq:Ah}) is
well-known to be of order~2, both in space and time.
Furthermore, this result holds as soon as the exact solution admits Taylor
approximations up to order~4.
Therefore, we formally prove Lemma~\ref{lm:consistency}, where sufficient
regularity is defined in Section~\ref{sec:o}, and uniform consistency is
defined by Equation~(\ref{eq:consistency}).

\begin{lemma}
  \label{lm:consistency}
  If the exact solution of the wave equation defined by Equations~(\ref{eq:L})
  to~(\ref{eq:A}) is sufficiently regular of order~4 uniformly on
  $\Omega\times[0,\tmax]$, then the numerical scheme defined by
  Equations~(\ref{eq:Lh}) to~(\ref{eq:Ah}) is consistent with the continuous
  problem at order (2,~2) uniformly on interval $[0,\tmax]$.
\end{lemma}

The proof uses properties of uniform Taylor approximations.
It deals with long and complex expressions, but it is still straightforward as
soon as discrete quantities are banned in favor of uniform continuous
quantities that encompass them.
Indeed, the key point is to use Taylor approximations that are uniform on a
compact set including all points of all grids showing up when the
discretization parameters tend to zero.

For instance, to handle the initialization phase of Equation~(\ref{eq:L1h}), we
consider a uniform Taylor approximation of order~1 of the following function
(for any~$v$ sufficiently regular of order~3):
\begin{multline*}
  ((x, t), (\Delta x, \Delta t)) \mapsto \mymultlinenewline
  \frac{v (x, t + \Delta t) -
    v (x, t)}{\Delta t} -
  \frac{\Delta t}{2} c^2 \frac{v (x + \Delta x, t) -
    2 v (x, t) +
    v (x - \Delta x, t)}{\Delta x^2}.
\end{multline*}
Thanks to the second-order approximation used in Equation~(\ref{eq:L1h}), the
initialization phase does not deteriorate the consistency order.

Note that, to provide Taylor approximations for both
expressions~$v(x+\Delta x,t)$ and~$v(x-\Delta x,t)$, the formal definition of
the uniform Taylor approximation via Equation~(\ref{eq:bigO}) must accept
negative changes in~$\bfx$, and not only the positive space grid steps, as a
naive specification would state.

\subsection{Stability}
\label{sec:stab}

When showing convergence of a numerical scheme, energy-based techniques use a
specific stability statement involving an estimation of the discrete energy.
Therefore, we formally prove Lemma~\ref{lm:energystability} where the
CFL$(\xi)$ condition is defined by Equations~(\ref{eq:cn}) and~(\ref{eq:cfl}).

\begin{lemma}
  \label{lm:energystability}
  For all $\xi\in]0,1[$, if discretization parameters satisfy the CFL$(\xi)$
  condition, then the numerical scheme defined by Equations~(\ref{eq:Lh})
  to~(\ref{eq:Ah}) is energetically stable uniformly on interval $[0,\tmax]$.
  Moreover, constants in Equation~(\ref{eq:energystability}) defining uniform
  energy stability are
  \begin{equation}
    \label{eq:energystabilityconstants}
    C_1 = \sqrt{E_{\rm h} (c) (p_{\rm h})^{\demi}} \qquad \mbox{and} \qquad
    C_2 = \frac{1}{\sqrt{2 \xi (2 - \xi)}}.
  \end{equation}
\end{lemma}

First, we show the variation of the discrete energy between two consecutive
time steps to be proportional to the source term.
In particular, the discrete energy is constant when the source is inactive.
Then, we establish the following lower bound for the discrete energy:
\begin{equation*}
  \demi (1 - \CN^2)
  \norme{\Delta x}{%
    \left(i \mapsto \frac{p_i^{k + 1} - p_i^k}{\Delta t}\right)}^2
  \le {E_{\rm h} (c) (p_{\rm h})^{k + \demi}}
\end{equation*}
where~$\CN$ is the Courant number defined by Equation~(\ref{eq:cn}), and
$k=k_{\Delta t}(t)$ for some~$t$ in~$[0,\tmax]$.
Therefore, the discrete energy is a positive semidefinite quadratic form, since
its nonnegativity directly follows the CFL$(\xi)$ condition.
Finally, Lemma~\ref{lm:energystability} follows in a straightforward manner
from the estimate:
\begin{equation*}
  \norme{\Delta x}{(i \mapsto p_i^{k + 1} - p_i^k)} \le
  2 C_2 \Delta t {E_{\rm h} (c) (p_{\rm h})^{k + \demi}}
\end{equation*}
where~$C_2$ is given in Equation~(\ref{eq:energystabilityconstants}).

Note that this stability result is valid for any discrete Cauchy
data~$u_{0,{\rm h}}$ and~$u_{1,{\rm h}}$, so it may be suboptimal for specific
choices.

\subsection{Convergence}
\label{sec:conv}

We formally prove Theorem~\ref{th:convergence}, where sufficient regularity is
defined in Section~\ref{sec:o}, the CFL$(\xi)$ condition is defined by
Equations~(\ref{eq:cn}) and~(\ref{eq:cfl}), and uniform convergence is defined
by Equation~(\ref{eq:convergence}).

\begin{theorem}
  \label{th:convergence}
  For all $\xi\in]0,1[$, if the exact solution of the wave equation defined by
  Equations~(\ref{eq:L}) to~(\ref{eq:A}) is sufficiently regular of order~4
  uniformly on $\Omega\times[0,\tmax]$, and if the discretization parameters
  satisfy the CFL($\xi$) condition, then the numerical scheme defined by
  Equations~(\ref{eq:Lh}) to~(\ref{eq:Ah}) is convergent of order (2,~2)
  uniformly on interval $[0,\tmax]$.
\end{theorem}

First, we prove that the convergence error~$e_{\rm h}$ is itself solution of
the same numerical scheme but with different input data.
(Of course, there is no associated continuous problem here.)
In particular, the source term on the right-hand side of Equation~(\ref{eq:Lh})
is here the truncation error~$\eps_{\rm h}$ associated with the numerical
scheme for~$p_{\rm h}$.
Then, from Lemma~\ref{lm:energystability} (stating uniform energy stability),
we have an estimation of the discrete energy associated with the convergence
error $E_{\rm h}(c)(e_{\rm h})$ that involves the sum of the corresponding source
terms, {\ie} the truncation error.
Finally, from Lemma~\ref{lm:consistency} (stating uniform consistency), this
sum is shown to have the same growth rate as $\Delta x^2+\Delta t^2$.

\section{Formal Proof of Boundedness of the Round-off Error}
\label{sec:round}

Beyond the proof of convergence of the numerical scheme, which is here
mechanically checked, we can now go into the sketch of the main steps of the
proof of the boundedness of the round-off error due to the limited precision of
computations in the program.
Full formal proofs are available from~\cite{code_and_proofs}.
Feedback about the novelty of this part is discussed in
Section~\ref{sec:mathsandproofs}.

For the {\clang}~program presented in Section~\ref{sec:wave}, naive
\gloss{forward error analysis} gives an error bound that is proportional to
$2^{k-53}$ for the computation of any~$p_i^k$.
If this bound was tight, it would cause the numerical scheme to compute only
noise after a few steps.
Fortunately, round-off errors actually cancel themselves in the program.
To take into account cancellations and hence prove a usable error bound, we
need a precise statement of the round-off error~\cite{Bol09b} to exhibit the
cancellations made by the numerical scheme.

The formal proof uses a comprehensive formalization of the IEEE-754 standard,
and hence both normal and subnormal ({\ie} very small) numbers are handled.
Other floating-point special values (infinities, NaNs) are proved in
Section~\ref{sec:auto} not to appear in the program.

\newcommand{\ul}{\underline}
\newcommand{\fl}{{\rm fl}}

\subsection{Local Round-off Errors}
\label{sec:local-round-off-error}

Let~$\delta_i^k$ be the local floating-point error that occurs during the
computation of~$p_i^k$.
For $k=0$ (resp. $k=1$, and $k\geq 2$), the local error corresponds to the
floating-point error at Line~9 of Listing~\ref{l:unannotated_code} (resp. at
Lines 19--20, and at Lines~32--33).
To distinguish them from the discrete values of previous sections, the
floating-point values as computed by the program are underlined below.
Quantities $\ul{a}$ and $\ul{p}_i^k$ match the expressions \codeacsl{a} and
\codeacsl{p[i][k]} in the \gloss{annotations}, while $a$ and $p_i^k$ are
represented by \lstinline[style=ACSL]{\exact(a)} and
\lstinline[style=ACSL]{\exact(p[i][k])}, as described in
Section~\ref{sec:tool_fp}.

The local floating-point error is defined as follows,
for all $k$ in $[1..\kmax-1]$, for all $i$ in $[1..\imax-1]$:
\begin{eqnarray*}
  \delta_i^{k+1} & = &
  \left(2 \ul{p}_i^k - \ul{p}_i^{k-1} +
    a (\ul{p}_{i+1}^k - 2 \ul{p}_i^k + \ul{p}_{i-1}^k)\right)
  - \ul{p}_i^{k+1} , \\
  \delta_i^1 & = &
  \left(\ul{p}_i^0 +
    \frac{a}{2} (\ul{p}_{i+1}^0 - 2 \ul{p}_i^0 + \ul{p}_{i-1}^0)\right)
  - \ul{p}_i^1
  - \left(\delta_i^0
    + \frac{a}{2} \left(\delta_{i+1}^0 - 2 \delta_i^0 + \delta_{i-1}^0\right)
  \right),\\
  \delta_i^0 & = & p_i^0 - \ul{p}_i^0.
\end{eqnarray*}
Note that the program presented in Section~\ref{sec:program} gives us that,
for all $k$ in $[1..\kmax-1]$, for all $i$ in $[1..\imax-1]$:
\begin{eqnarray*}
  \ul{p}_i^{k+1} & = &
  \fl\left(2 \ul{p}_i^k - \ul{p}_i^{k-1} +
    a (\ul{p}_{i+1}^k - 2 \ul{p}_i^k + \ul{p}_{i-1}^k)\right), \\
  \ul{p}_i^1 & = &
  \fl\left(\ul{p}_i^0 +
    \frac{a}{2} (\ul{p}_{i+1}^0 - 2 \ul{p}_i^0 + \ul{p}_{i-1}^0)\right), \\
  \ul{p}_i^0 & = &
  \fl\left(p_0(i \ul{\Delta x})\right).
\end{eqnarray*}
where $\fl(\cdot)$ means that all the arithmetic operations that appear between
the parentheses are actually performed by floating-point arithmetic, hence a
bit off.

In order to get a bound on $\delta_i^k$, we need to know the range of
$\ul{p}_i^k$.
For the bound to be usable in our correctness proof, the range has to be
$[-2,2]$.
We assume here that the exact solution is bounded by~$1$;
if not, we normalize by the maximum value and use the linearity of the problem.
We have proved this range by induction on a simple triangle inequality taking
advantage of the fact that at each point of the grid, the floating-point value
computed by the program is the sum of the exact solution, the method error, and
the round-off error.

To prove the bound on $\delta_i^k$, we perform forward error analysis and then
use \gloss{interval arithmetic} to bound each intermediate
error~\cite{DauMel10}.
We prove that, for all $i$ and $k$, we have
$|\delta_i^{k}| \le 78 \cdot 2^{-52}$ for a reasonable error bound for $a$,
that is to say $|\ul{a} - a| \le 2^{-49}$.

\subsection{Global Round-off Error}
\label{sec:global-round-off-error}

Let $\Delta_i^k=\ul{p}_i^k-p_i^k$ be the global floating-point error on~$p_i^k$.
The global floating-point error depends not only on the local floating-point
error at position $(i,k)$, but also on all the local floating-point errors
inside the space-time dependency cone of apex $(i,k)$.

From the linearity of the numerical scheme, the global floating-point error is
itself solution of the same numerical scheme with discrete inputs corresponding
to the local floating-point error:
\[u_{0,{\rm h}}=\delta_{\rm h}^0=0,\quad
u_{1,{\rm h}}=\frac{\delta_{\rm h}^1}{\Delta t},\quad
s_{\rm h}=\left(k\mapsto\frac{\delta_{\rm h}^{k+1}}{\Delta t^2}\right).\]

Taking advantage of the remark about image theory made in
Section~\ref{sec:discrete}, we see the global floating-point error as the
restriction to the domain~$\Omega$ of the solution of the same numerical
scheme set on the entire real axis without the Dirichlet boundary condition of
Equation~(\ref{eq:dirh}), and with extended discrete inputs.
Therefore, the expression of the global floating-point error is given by the
convolution of the discrete fundamental solution on the entire real axis and
the extended local floating-point error.

Let us denote by~$\lambda_{\rm h}$ the discrete fundamental solution.
It is solution of the same numerical scheme set on the entire real axis with
null discrete inputs except $u_{1,0}=\frac{1}{\Delta t}$.
Then, we can state the following result. (See~\cite{Bol09b} for a direct proof
that does not follow above remarks.)

\begin{theorem}
\label{th:delta}
\[\forall k \geq 0,\, \forall i \in [0..\imax], \quad
\Delta_i^k = \sum_{l=0}^k \sum_{j=-l}^l
\tilde{\delta}_{i-j}^{k-l} \, \lambda_j^{l+1}.\]
\end{theorem}

Note that, for all $k$, we have $\Delta_0^k=\Delta_{\imax}^k=0$.
Note also that $\lambda_i^k$ vanishes as soon as $|i|\geq k$, hence the sum
over~$j$ could be replaced by an infinite sum over all integers.

\subsection{Bound on the Global Round-off Error}

The analytic expression of $\Delta_i^k$ can be used to obtain a bound on the
round-off error.
We will need two lemmas for this purpose.

\begin{lemma} \label{lm:l1}
\[\forall k \geq 0, \quad
\sigma^k\eqdef\sum_{i=-\infty}^{+\infty} \lambda_i^k=k.\]
\end{lemma}

\begin{myproof}
The sum satisfies the following linear recurrence:
for all $k\geq 1$, $\sigma^{k+1}-2\sigma^k+\sigma^{k-1}=0$.
Since $\sigma^0=0$, and $\sigma^1=1$, we have, for all $k\geq 0$,
$\sigma^k=k$.
\end{myproof}

\begin{lemma}
\label{th:jacobi}
\[\forall k \geq 0,\, \forall i, \quad
\lambda_i^k \ge 0.\]
\end{lemma}

The sketch of the proof is given in \myappendixreftext\ref{sec:jacobi}.
It uses complex algebraic results (the positivity of sums of Jacobi
polynomials), and was not formalized in {\coq}.

\begin{theorem}
\label{th:roundoff_error}
\[\forall k \geq 0,\, \forall i \in [0..\imax], \quad
\left| \Delta_i^k \right| \le 78 \cdot 2^{-53} (k+1) (k+2).\]
\end{theorem}

\begin{myproof}
According to Theorem~\ref{th:delta}, $\Delta_i^k$ is equal to
$\sum_{l=0}^k \sum_{j=-l}^l \lambda_j^{l+1} \ \tilde{\delta}_{i-j}^{k-l}$.
We know from Section~\ref{sec:local-round-off-error} that for all $j$ and $l$,
$|\tilde{\delta}_j^l| \le 78 \cdot 2^{-52}$, and from Lemma~\ref{lm:l1} that
$\sum \lambda_i^{l+1}=l+1$.
Since the $\lambda$'s are nonnegative (Lemma~\ref{th:jacobi}), the error is
easily bounded by $78 \cdot 2^{-52} \sum_{l=0}^k (l+1)$.
\end{myproof}

Except for Lemma~\ref{th:jacobi}, all proofs described in this section have
been machine-checked using {\coq}.
In particular, the proof of the bound on $\delta_i^k$ was done automatically by
calling Gappa from {\coq}.
Lemma~\ref{th:jacobi} is a technical detail compared to the rest of our work,
and is not worth the tremendous work it would require to be implemented in
{\coq}: keen results on integrals but also definitions and results about the
Legendre, Laguerre, Chebychev, and Jacobi polynomials (see
\myappendixreftext\ref{sec:jacobi}).

\section{Formal Proof of the {\clang}~Program}
\label{sec:program_proof}

To completely prove the actual {\clang}~program, we bound the total error
incorporating both the method error and the round-off error.
Then, we show how we prove the absence of coding errors (such as overflow and
out-of-bound access) and how we relate the behavior of the program to the
preceding {\coq} proofs.

\subsection{Total Error}
\label{sec:total_error}

Let~$\cal E_{\rm h}$ be the total error.
It is the sum of the method error (or convergence error) $e_{\rm h}$ of
Sections~\ref{sec:wave_conv} and~\ref{sec:conv}, and of the round-off
error~$\Delta_h$ of Section~\ref{sec:round}.

From Theorem~\ref{th:roundoff_error}, we can estimate%
\footnote{When $\tau\geq 2$, we have $(\tau+1)(\tau+2)\leq 3\tau^2$.}
the following upper bound for the spatial norm of the round-off error when
$\Delta x\leq 1$ and $\Delta t\leq\tmax/2$:
for all $t\in[0,\tmax]$,
\begin{eqnarray*}
  \ndx{\left(i \mapsto \Delta_i^{k_{\Delta t}(t)}\right)}
  & = & \sqrt{\sum_{i=0}^{\imax}
    \left( \Delta_i^{k_{\Delta t}(t)} \right)^2 \Delta x} \\
  & \le & \sqrt{(\imax+1)\Delta x} \cdot 78  \cdot 2^{-53}
  \left( \frac{\tmax}{\Delta t} + 1 \right)
  \left( \frac{\tmax}{\Delta t} + 2 \right) \\
  & \leq & \sqrt{\xmax - \xmin + 1} \cdot 78 \cdot 2^{-53} \cdot
  3 \frac{\tmax^2}{\Delta t^2}.
\end{eqnarray*}

Thus, from the triangle inequality for the spatial norm, we obtain the
following estimation of the total error:
\begin{multline}
  \forall t \in [0, \tmax], \;
  \forall \bfdeltax, \quad
  \|\bfdeltax\| \leq \min(\alpha_e, \alpha_\Delta) \Rightarrow \\
  \label{eq:totalerrorbound}
  \ndx{\left(i \mapsto {\cal E}_i^{k_{\Delta t}(t)}\right)} \leq
  C_e (\Delta x^2 + \Delta t^2) + \frac{C_\Delta}{\Delta t^2}
\end{multline}
where the method error constants~$\alpha_e$ and~$C_e$ were extracted from the
{\coq} proof (see Section~\ref{sec:conv}) and are given in terms of the
constants for the Taylor approximation of the exact solution at degree~3
($\alpha_3$ and $C_3$), and at degree~4 ($\alpha_4$ and $C_4$) by
\begin{eqnarray}
  \label{eq:alpha_e}
  \alpha_e & = & \min (1, \tmax, \alpha_3, \alpha_4), \\
  \label{eq:C_e}
  C_e & = & 4 C_2 \tmax \sqrt{\xmax - \xmin}
  \left(\frac{C^\prime}{\sqrt{2}} + 2 C_2 (\tmax + 1) C^{\prime\prime} \right)
\end{eqnarray}
with~$C_2$ coming from Equation~(\ref{eq:energystabilityconstants}),
$C^\prime=\max(1,C_3+c^2 C_4+1)$,
$C^{\prime\prime}=\max(C^\prime,2(1+c^2)C_4)$,
and where the round-off constants~$\alpha_\Delta$ and~$C_\Delta$, as explained
above,  are given by
\begin{eqnarray}
  \label{eq:alpha_Delta}
  \alpha_\Delta & = & \min(1, \tmax/2), \\
  \label{eq:C_Delta}
  C_\Delta & = & 234 \cdot 2^{-53} \tmax^2 \sqrt{\xmax - \xmin + 1}.
\end{eqnarray}

Of course, decreasing the size of the grid step decreases the method error, but
at the same time, it increases the round-off error.
Therefore, there exists a minimum for the upper bound on the total error,
corresponding to optimal grid step sizes that may be determined using above
formulas.

On specific examples, one observes that this upper bound on the total error can
be highly overestimated (see Section~\ref{sec:mathsandproofs}).
Nevertheless, one also observes that the asymptotic behavior of the upper bound
of the total error for high values of the grid steps is close to the asymptotic
behavior of the effective total error~\cite{BCFMMW12}.

\subsection{Automation and Manual Proofs}
\label{sec:auto}

\begin{figure}[htb]
\begin{center}
  \includegraphics[width=\linewidth]{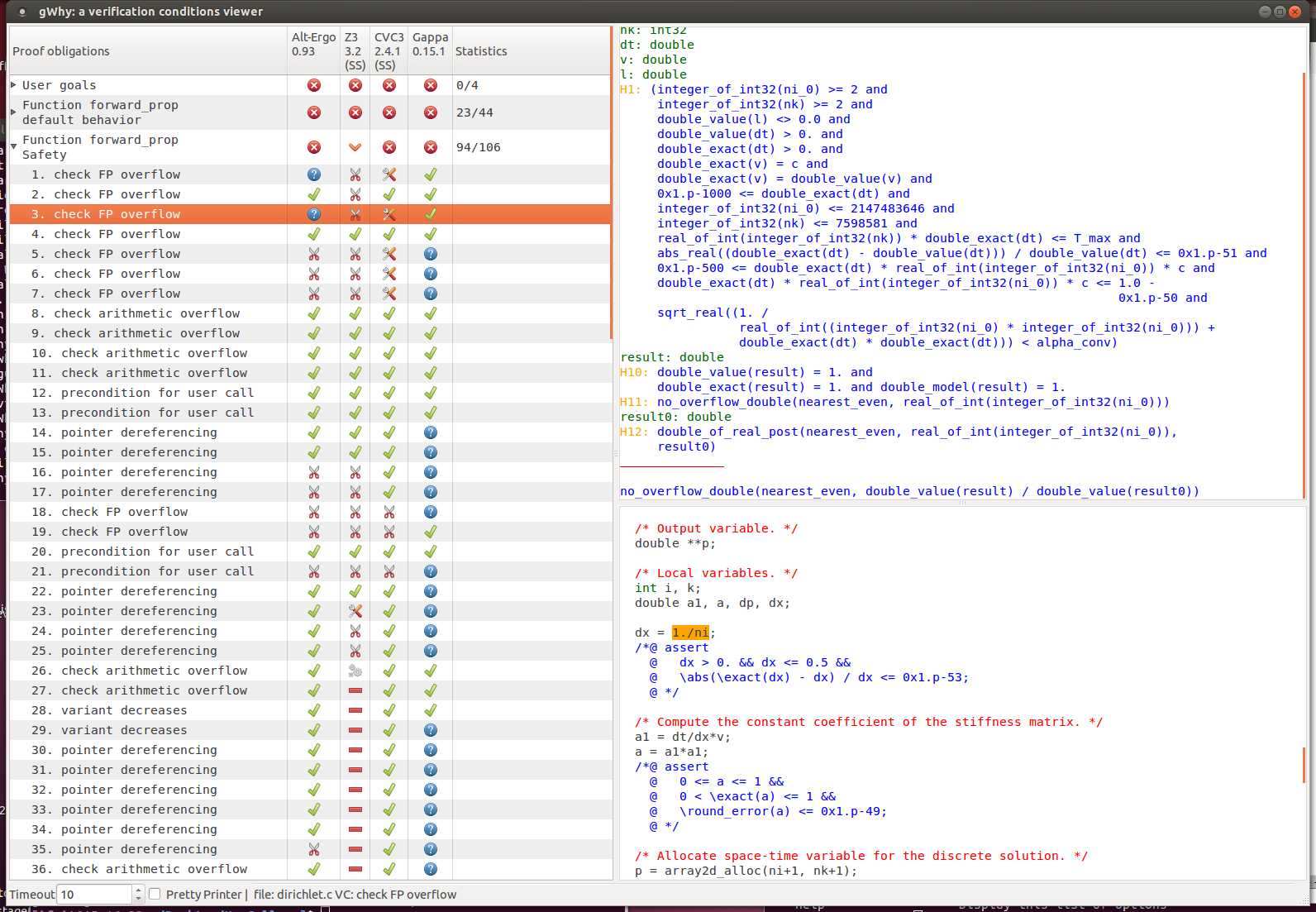}
\end{center}
\caption{Screenshot of gWhy.
  Left: list of all verification conditions (VCs) and their proof status
  with respect to the various automatic tools.
  Upper right: statement of the currently selected VC.
  Lower right: location in the source code where this VC originates from.}
\label{fig:screenshot}
\end{figure}

Given the program code, the Frama-C/Jessie/Why tools generate 150
\gloss{verification conditions} that have to be proved (see
Section~\ref{sec:tool_fp}).
While possible, proving all of them in {\coq} would be rather tedious.
Moreover, systematically using {\coq} would lead to a rather fragile construct:
any later modification to the program, however small it is, would cause
different proof obligations to be generated, which would then require
additional human interaction to adapt the {\coq} proofs.
We prefer to have \gloss{automated theorem provers} (\gloss{SMT solvers} and
Gappa, see Sections~\ref{sec:tool_fp} and~\ref{sec:tool_gappa}) prove as many
of them as possible, so that only the most intricate ones are left to be proved
in {\coq}.

Figure~\ref{fig:screenshot} displays a screenshot of the gWhy graphical
interface.
For instance, any discrepancy between the code and its \gloss{specification}
would be highlighted on the lower right part when the corresponding
verification condition is selected.
A~comprehensive view on all the verification conditions and how they are proved
can be found in~\cite{code_and_proofs}.

Verification conditions split into safety \gloss{goals} and behavior goals.
Safety goals are always generated, even in the absence of any specification.
Proving them ensures that the program always successfully terminates.
They check that matrix accesses are in range, that the loop variants decrease
and are nonnegative (thus loops terminate), that integer and floating-point
arithmetic operations do not overflow, and so on.
On such goals, automatic provers are helpful, as they prove about 90\,\% of
them.

Behavior goals relate the program to its specification.
Proving them ensures that, if the program terminates without error, then it
returns the specified result.
They check that \gloss{loop invariants} are preserved, that assertions hold,
that preconditions hold before function calls, that postconditions are implied
by preconditions, and so on.
Automatic provers are a bit less helpful, as they fail for half of the behavior
goals.
That is why we resort to an \gloss{interactive higher-order theorem prover},
namely {\coq}.

{\coq} developments split into two sets of files.
The first set, for a total of more than 4,000 lines of script, is dedicated to
the proof of the boundedness of the method error, {\ie} the convergence of the
numerical scheme.
About half of it is generic, meaning it can be re-used for other PDEs or ODEs:
definitions and lemmas about big~O, scalar product, $\R^n$ and so on.
The other half is dedicated to the 1D wave equation and the chosen numerical
scheme.
The second set, for a total of more than 12,000 lines, deals with both the
proof of the boundedness of the round-off error, and the proof of the absence
of runtime errors.
Note that slightly less than half of those {\coq} developments are statements
of the second set that are automatically generated by the Frama-C/Jessie/Why
framework, whereas the other half---statements for the convergence and all
proofs---is manually edited.

\section{Discussion}
\label{sec:discussion}

We emphasize now the difficulties and the originality of this work that may
strike the computational scientist's mind when confronted to the formal proof
of a numerical program.
Finally, we discuss some future work.

\subsection{Applied Mathematicians and Formal Proof of Programs}
\label{sec:mathsandproofs}

From the point of view of computational scientists, the big surprise of a
\gloss{formal proof} development may come from the requirement of writing an
extremely detailed conventional pen-and-paper proof.
Indeed, an \gloss{interactive proof assistant} such as {\coq} is not able to
elaborate a sketch of proof, nor to automatically find a proof of a lemma---but
for the most trivial facts.
Therefore, a formal proof is completely human driven, up to the utmost detail.
Fortunately, we do not need to specify all the mathematical facts from scratch,
since {\coq} provides a collection of known facts and theories that can be
reused to establish new results.
However, classical results from Hilbert space analysis, uniform asymptotic
comparison of functions, and Taylor approximations have not yet been proved in
{\coq}.
Thus, we had to develop them, as well as in the detailed pen-and-paper proof.
In the end, the detailed pen-and-paper proof is about 50 pages long.
Surprisingly, it is roughly of the same size as the complete {\coq} formal
proof: both collections of source files are approximately 5~kLOC long, and
weigh about 160~kB.

The case of the big~O notation deserves some comments.
Typically, this is a notion about which pen-and-paper proofs do not elaborate
on.
Of course, misinterpreting the involved existential constants leads to
erroneous reasoning.
This issue is seldom dealt with in the literature; however, we point out an
informal mention in~\cite{tho:npd:95} with a relevant remark about pointwise
consistency versus norm consistency. 
We had to focus on the uniform notion of big~O equalities in~\cite{BCFMMW10} in
the context of the infinite string, for which the space domain~$\Omega$ is the
whole real axis.
In the present paper, we deal with finite strings.
Thus, for compactness reasons, both uniform and nonuniform notions of big~O
notations are obviously equivalent.
Nevertheless, we still use the more general uniform big~O notion to share most
of the proof developments between finite and infinite cases.

With pen-and-paper proofs, it is uncommon to expand constants involved in the
estimation of the method error (namely~$\alpha$, and~$C$ in the big~O equality
for the norm of the convergence error).
In contrast, {\coq} can automatically extract from the formal proof the
constants appearing in the upper bound for the total error, which also takes
the round-off error into account.
The mathematical expressions given in Section~\ref{sec:total_error} are as
sharp as the estimations used in the proof can provide.
The same expressions could be obtained by hand with extreme care from the
pen-and-paper proof, at least for the method error contribution.
However, in specific examples, the actual total error may be much smaller than
the established estimation~\cite{BCFMMW12}.
This is essentially due to the use of {\apriori} error estimations, worsened by
the energy-based techniques which are known to become less precise as the
CFL condition becomes optimal ({\ie} as the Courant number goes up to~1).
Indeed, constant~$C_e$ in Equation~(\ref{eq:C_e}) grows as~$1/\xi$ for small
values of~$\xi$.
An alternative to obtain more accurate {\apriori} bounds is to use a leapfrog
scheme for the equivalent first-order system~\cite{bec:esn:09}.

Nowadays, {\aposteriori} error estimates have become popular tools to control
the method error, {\eg} see~\cite{ver:rap:96,or:vvs:10}.
Indeed, the measure at runtime of the error due to the different stages of
approximation provides much tighter bounds than any {\apriori} technique. 
For instance, the approach allows to adapt stopping criteria of internal
iterative solvers to the overall numerical error, and to save a significant
amount of computations~\cite{ev:ain:14}.
From the certification point of view, such developments present themselves as
part of the numerical method, as a complement to the numerical scheme:
{\aposteriori} error estimates are first designed and proved on paper, then
implemented into the program, and finally make use of floating-point operators
in their computations.
They deserve exactly the same care in terms of code \gloss{verification}:
first, their definition and behavior must be specified in the
\gloss{annotations}, and then the generated \gloss{verification conditions}
must be formally proved.

\bigskip

The other big surprise for computational scientists is the possibility to take
into account the round-off error cancellations occurring during the
computations, and to evaluate a sensible bound that encompasses all aspects of
the IEEE-754 floating-point arithmetic.
Of course, this is typically the kind of results that are not accessible
outside the logical framework of a mechanical interactive proof assistant as
{\coq}.

Furthermore, in numerical analysis, one usually evaluates bounds on the
absolute value of quantities of interest.
To obtain an upper bound on the round-off error, we needed a result about the
sign of the fundamental solution of the discrete scheme, not about a bound for
its absolute value.
To our best knowledge, the only way to capture such a result is purely
algebraic:
the closed-form expression of the solution is first obtained using a generating
function, then it may be recognized as a combination of Jacobi polynomials that
happens to be nonnegative (see \myappendixreftext\ref{sec:jacobi}).

Such a result about the sign of the discrete fundamental solution may not hold
for other numerical schemes.
In our case, just assuming a bound on the absolute value of the discrete
fundamental solution would lead to an estimation of the global round-off error
about 4~times higher than that of Theorem~\ref{th:roundoff_error}.

\bigskip

Finally, differential equations introduce an issue in the annotations.
Due to the underlying logic, the annotations have to define the solution of the
PDE by using \gloss{first-order formulas} stating differentiability, instead of
second-order formulas involving differential operators.
This makes the annotations especially tedious and verbose.
In the end, the {\clang} code grows by a factor of about~4, as can be seen when
comparing the original source in Listing~\ref{l:unannotated_code} with the
fully annotated source in \myappendixreftext\ref{sec:annotated_source}.

\subsection{Future Work}
\label{sec:future_work}

Given its cost, we do not plan to apply the same methodology to any other
scheme or any other PDE at random.
We nevertheless want to take advantage of our acquired expertise on several
related topics.

In floating-point arithmetic, an algorithm is said to be
{\em numerically stable} if it has a small \gloss{forward error}.
This property is certainly related to the notion of stability used in
scientific computing for numerical schemes, which states that computed values
do not diverge.
As a rule of thumb, stable schemes are considered numerically stable.
We plan to formally prove this statement for a large class of numerical
schemes.
This will require the formalization of the notion of numerical scheme, and of
the property of stability.
Indeed, the definition given in Section~\ref{sec:wave_conv} is dedicated to
numerical schemes for second-order in time evolution problems.

The mechanism used to extract the constants involved in the estimation of the
total error can be deployed on a wider scale to generate the whole program.
Instead of formally proving properties specified in an existing program, one
formalizes a problem, proves it has a solution, and obtains for free an
{\ocaml} program by automatic extraction from the {\coq} proof.
Such a program is usually inefficient, but it is zero-defect, provided it is
not modified afterward.
For computational science problems involving real numbers, efficiency is not the
only issue: the extracted program would also be hindered by the omnipresence of
classical real numbers in the formalization.
Therefore, being able to extract efficient and usable zero-defect programs
would be an interesting long-term goal for critical numerical problems.

In the present exploratory work, we consider the simple second-order centered
finite difference scheme for the one-dimensional wave equation.
Further works involve scaling to higher-order numerical schemes, and/or to
higher-dimensional problems.
Finite support functions, and summations on such support played a much more
important role in the {\coq} developments than we initially expected.
Therefore, we consider using the {\ssreflect} interface and libraries for
{\coq}~\cite{BGOP08}, so as to simplify manipulations of these objects in
higher-dimensional cases.

Another long-term perspective is to generalize the mathematical notions
formalized in {\coq} to be able to apply the approach to finite difference
schemes for other PDEs, and ODEs.
Steps in this direction concern the generalization of the use of convolution
with the discrete fundamental solution to bound the round-off error, and the
formal proof of the Lax equivalence theorem.

Finally, a more ambitious perspective is to formally deal with the finite
element method.
This first requires a {\coq} formalization of mathematical tools from Hilbert
space analysis on Sobolev spaces such as the Lax-Milgram theorem (for existence
and uniqueness of exact and approximated solutions), and the handling of meshes
(that may be regular or not, and structured or not).
Then, the finite element method may be proved convergent in order to formally
guarantee a large class of numerical analysis programs that solve linear PDEs
on complex geometries.
The purpose is to go beyond Laplace's equation set on the unit square, {\eg}
handle mixed boundary conditions, and extend to mixed and mixed hybrid finite
element methods.

\section{Conclusion}
\label{sec:conclusion}

We have presented a comprehensive formal proof of a {\clang}~program solving
the 1D acoustic wave equation using the second-order centered finite difference
explicit scheme. Our proof includes the following aspects.
\begin{itemize}
\item \emph{Safety}: we prove that the {\clang}~program terminates and is free
  of runtime failure such as division by zero, array access out of bounds, null
  pointer dereference, or arithmetic overflow.
  The latter includes both integer and floating-point overflows.
\item \emph{Method error}: we show that the numerical scheme is convergent of
  order $(2,2)$ uniformly on the time interval, under the assumptions that the
  exact solution is sufficiently regular ({\ie} its Taylor polynomials of order
  up to~4 are uniform approximations on the space--time domain) and that the
  Courant-Friedrichs-Lewy condition holds.
\item \emph{Round-off error}: we bound all round-off errors resulting from
  floating-point computations used in the program.
  In particular, we show how some round-off errors cancel, to eventually get a
  meaningful bound.
\item \emph{Total error}: altogether, we are able to provide explicit (and
  formally proved) bounds for the sum of method and round-off errors.
\end{itemize}
To our knowledge, this is the first time such a comprehensive proof is
achieved.
Related to the small length of the program (a few tens of lines), the total
cost of the formal proof is huge, if not frightening: several man-months of
work, three times more annotations to be inserted in the program than lines of
code, over 16,000 lines of {\coq} proof scripts, 30 minutes of CPU time to
check them on a 3-GHz processor.

The recent introduction of formal methods in DO-178C%
\footnote{%
  DO-178C is the newest version of Software Considerations in Airborne Systems
  and Equipment Certification, which is used by national certification
  authorities.
}
shows the need for the verification of numerical programs in the context of
embedded critical software.
Considering the work we have presented, one can hardly think of verifying
numerical codes on the scale of a large airborne system.
Yet we think our techniques and a large subset of our proof can be reused and
would significantly decrease the workload of such a proof.
This is to be combined with increased proof automation, so that user
interaction is minimized.
Finally, the challenge is to provide tools that are usable by computational
scientists that are not specialists of formal methods.

In conclusion, combining scientific computing and formal proofs is now
considered an important matter by logicians.
Formal tools for scientific computing are being actively developed and progress
is done to get them mature and usable for non specialists.
It seems that it is now time for computational scientists to take a keen
interest in this area.

\section*{Acknowledgments}
We are grateful to Manuel Kauers, Veronika Pillwein, and Bruno Salvy, who
provided us help with the nonnegativity of the fundamental solution of the
discrete wave equation (Lemma~\ref{th:jacobi}).
We are also thankful to Vincent Martin for his constructive remarks on this
article.
Last, we feel indebted to the reviewers for their invaluable suggestions.

\bibliography{biblio}

\begin{thebibliography}{10}

\bibitem{AnAsRo99}
George~E. Andrews, Richard Askey, and Ranjan Roy.
\newblock {\em Special functions}.
\newblock Cambridge University Press, Cambridge, 1999.

\bibitem{AsGa72}
Richard Askey and George Gasper.
\newblock Certain rational functions whose power series have positive
  coefficients.
\newblock {\em The American Mathematical Monthly}, 79:327--341, 1972.

\bibitem{CVC3}
Clark Barrett and Cesare Tinelli.
\newblock {CVC3}.
\newblock In {\em 19th International Conference on Computer Aided Verification
  (CAV '07)}, volume 4590 of {\em LNCS}, pages 298--302. Springer-Verlag, 2007.
\newblock {B}erlin, Germany.

\bibitem{ACSL}
Patrick Baudin, Pascal Cuoq, Jean-Christophe Filli\^atre, Claude March\'e,
  Benjamin Monate, Yannick Moy, and Virgile Prevosto.
\newblock {\em ACSL: ANSI/ISO C Specification Language, version 1.5}, 2009.
\newblock \url{http://frama-c.com/acsl.html}.

\bibitem{bec:esn:09}
\'Eliane B\'ecache.
\newblock {\'E}tude de sch\'emas num\'eriques pour la r\'esolution de
  l'\'equation des ondes.
\newblock Master Mod\'elisation et simulation, Cours ENSTA, 2009.
\newblock
  \url{http://www.ensta-paristech.fr/~becache/COURS-ONDES/Poly-num-0209.pdf}.

\bibitem{Coq}
Yves Bertot and Pierre Cast\'eran.
\newblock {\em Interactive Theorem Proving and Program Development. {Coq'Art}:
  The Calculus of Inductive Constructions}.
\newblock Texts in Theoretical Computer Science. Springer, 2004.

\bibitem{BGOP08}
Yves Bertot, Georges Gonthier, Sidi Ould~Biha, and Ioana Pasca.
\newblock Canonical big operators.
\newblock In {\em 21st International Conference on Theorem Proving in Higher
  Order Logics (TPHOLs'08)}, volume 5170 of {\em LNCS}, pages 86--101.
  Springer, 2008.

\bibitem{Bol04b}
Sylvie Boldo.
\newblock {\em Preuves formelles en arithm\'etiques \`a virgule flottante}.
\newblock PhD thesis, \'Ecole Normale Sup\'erieure de Lyon, 2004.

\bibitem{Bol09b}
Sylvie Boldo.
\newblock Floats \& {R}opes: a case study for formal numerical program
  verification.
\newblock In {\em 36th International Colloquium on Automata, Languages and
  Programming}, volume 5556 of {\em LNCS - ARCoSS}, pages 91--102. Springer,
  2009.

\bibitem{BCFMMW10}
Sylvie Boldo, Fran\c{c}ois Cl{\'e}ment, Jean-Christophe Filli{\^a}tre, Micaela
  Mayero, Guillaume Melquiond, and Pierre Weis.
\newblock Formal proof of a wave equation resolution scheme: the method error.
\newblock In Matt Kaufmann and Lawrence~C. Paulson, editors, {\em 1st
  Conference on Interactive Theorem Proving Conference (ITP 2010)}, volume 6172
  of {\em LNCS}, pages 147--162. Springer, 2010.

\bibitem{code_and_proofs}
Sylvie Boldo, Fran\c{c}ois Cl{\'e}ment, Jean-Christophe Filli{\^a}tre, Micaela
  Mayero, Guillaume Melquiond, and Pierre Weis.
\newblock Annotated source code and {C}oq proofs, 2012.
\newblock \url{http://fost.saclay.inria.fr/wave_total_error.html}.

\bibitem{BCFMMW12}
Sylvie Boldo, Fran\c{c}ois Cl{\'e}ment, Jean-Christophe Filli{\^a}tre, Micaela
  Mayero, Guillaume Melquiond, and Pierre Weis.
\newblock Wave equation numerical resolution: a comprehensive mechanized proof
  of a {C} program.
\newblock {\em Journal of Automated Reasoning}, 50(4):423--456, 2013.

\bibitem{BoldoFilliatre07}
Sylvie Boldo and Jean-Christophe Filli\^atre.
\newblock Formal verification of floating-point programs.
\newblock In {\em 18th IEEE International Symposium on Computer Arithmetic},
  pages 187--194, 2007.

\bibitem{BolFilMel09}
Sylvie Boldo, Jean-Christophe Filli{\^a}tre, and Guillaume Melquiond.
\newblock Combining {Coq} and {Gappa} for certifying floating-point programs.
\newblock In Jacques Carette, Lucas Dixon, Claudio~Sarcedoti Coen, and
  Stephen~M. Watt, editors, {\em 16th Calculemus Symposium}, volume 5625 of
  {\em LNAI}, pages 59--74, 2009.

\bibitem{BJLM13}
Sylvie Boldo, Jacques-Henri Jourdan, Xavier Leroy, and Guillaume Melquiond.
\newblock A formally-verified {C} compiler supporting floating-point
  arithmetic.
\newblock In {\em Proceedings of the 21th IEEE Symposium on Computer
  Arithmetic}, pages 107--115, 2013.

\bibitem{BLM12}
Sylvie Boldo, Catherine Lelay, and Guillaume Melquiond.
\newblock Improving real analysis in {C}oq: A user-friendly approach to
  integrals and derivatives.
\newblock In C.~Hawblitzel and D.~Miller, editors, {\em Proceedings of the
  Second International Conference on Certified Programs and Proofs}, number
  7679 in LNCS, pages 289--304. Springer, 2012.

\bibitem{BolMel11}
Sylvie Boldo and Guillaume Melquiond.
\newblock {Flocq}: A unified library for proving floating-point algorithms in
  {Coq}.
\newblock In Elisardo Antelo, David Hough, and Paolo Ienne, editors, {\em 20th
  IEEE Symposium on Computer Arithmetic}, pages 243--252, 2011.

\bibitem{conchon08entcs}
Sylvain Conchon, \'Evelyne Contejean, Johannes Kanig, and St\'ephane Lescuyer.
\newblock {CC(X)}: Semantical combination of congruence closure with solvable
  theories.
\newblock In {\em Post-proceedings of the 5th International Workshop on
  Satisfiability Modulo Theories ({SMT 2007})}, volume 198-2 of {\em Electronic
  Notes in Computer Science}, pages 51--69. Elsevier Science Publishers, 2008.

\bibitem{CIC}
Thierry Coquand and Christine Paulin-Mohring.
\newblock Inductively defined types.
\newblock In P.~Martin-Löf and G.Mints, editors, {\em Colog’88}, volume 417
  of {\em LNCS}. Springer-Verlag, 1990.

\bibitem{cfl:pde:67}
Richard Courant, Kurt Friedrichs, and Hans Lewy.
\newblock On the partial difference equations of mathematical physics.
\newblock {\em IBM Journal of Research and Development}, 11(2):215--234, 1967.

\bibitem{DauMel10}
Marc Daumas and Guillaume Melquiond.
\newblock Certification of bounds on expressions involving rounded operators.
\newblock {\em Transactions on Mathematical Software}, 37(1):1--20, 2010.

\bibitem{DauRidThe01}
Marc Daumas, Laurence Rideau, and Laurent Th{\'e}ry.
\newblock A generic library for floating-point numbers and its application to
  exact computing.
\newblock In {\em Proceedings of the 14th International Conference on Theorem
  Proving in Higher-Order Logics (TPHOL'01)}, pages 169--184. Springer-Verlag,
  2001.

\bibitem{DinLauMel11}
Florent de~Dinechin, Christoph Lauter, and Guillaume Melquiond.
\newblock Certifying the floating-point implementation of an elementary
  function using {Gappa}.
\newblock {\em IEEE Transactions on Computers}, 60(2):242--253, 2011.

\bibitem{Z3}
Leonardo de~Moura and Nikolaj Bj{\o}rner.
\newblock {Z3}, an efficient {SMT} solver.
\newblock In {\em Proceedings of the 14th International Conference on Tools and
  Algorithms for the Construction and Analysis of Systems (TACAS'08)}, volume
  4963 of {\em Lecture Notes in Computer Science}, pages 337--340. Springer,
  2008.

\bibitem{Dij75}
Edsger~W. Dijkstra.
\newblock Guarded commands, nondeterminacy and formal derivation of programs.
\newblock {\em Communications of the ACM}, 18(8):453--457, 1975.

\bibitem{ev:ain:14}
Alexandre Ern and Martin Vohral\'ik.
\newblock Adaptive inexact newton methods with a posteriori stopping criteria
  for nonlinear diffusion pdes.
\newblock {\em SIAM Journal on Scientific Computing}, 35(4):A1761--A1791, 2014.

\bibitem{Fejer08}
L\'eopold Fej\'er.
\newblock Sur le d\'eveloppement d'une fonction arbitraire suivant les
  fonctions de {L}aplace.
\newblock {\em Comptes-Rendus de l'Acad\'emie des Sciences}, 146:224--227,
  1908.

\bibitem{filliatre07cav}
Jean-Christophe Filli\^atre and Claude March\'e.
\newblock The {Why/Krakatoa/Caduceus} platform for deductive program
  verification.
\newblock In {\em 19th International Conference on Computer Aided
  Verification}, volume 4590 of {\em LNCS}, pages 173--177. Springer, 2007.

\bibitem{Gasper77}
George Gasper.
\newblock Positive sums of the classical orthogonal polynomials.
\newblock {\em SIAM Journal on Mathematical Analysis}, 8(3):423--447, 1977.

\bibitem{Gold91}
David Goldberg.
\newblock What every computer scientist should know about floating-point
  arithmetic.
\newblock {\em ACM Comput. Surv.}, 23(1):5--48, 1991.

\bibitem{joh:pde:86}
Fritz John.
\newblock {\em Partial Differential Equations}.
\newblock Springer, 1986.

\bibitem{seL410}
Gerwin Klein, June Andronick, Kevin Elphinstone, Gernot Heiser, David Cock,
  Philip Derrin, Dhammika Elkaduwe, Kai Engelhardt, Rafal Kolanski, Michael
  Norrish, Thomas Sewell, Harvey Tuch, and Simon Winwood.
\newblock {seL4}: Formal verification of an operating system kernel.
\newblock {\em Communications of the ACM}, 53(6):107--115, 2010.

\bibitem{dal:rcf:47}
Jean {le Rond D'Alembert}.
\newblock Recherches sur la courbe que forme une corde tendue mise en
  vibrations.
\newblock In {\em Histoire de l'Acad\'emie Royale des Sciences et Belles
  Lettres (Ann\'ee 1747)}, volume~3, pages 214--249. Haude et Spener, Berlin,
  1749.

\bibitem{Ler09}
Xavier Leroy.
\newblock Formal verification of a realistic compiler.
\newblock {\em Communications of the ACM}, 52(7):107--115, 2009.

\bibitem{marche07plpv}
Claude March\'e.
\newblock Jessie: an intermediate language for {Java} and {C} verification.
\newblock In {\em Programming Languages meets Program Verification (PLPV)},
  pages 1--2. ACM, 2007.

\bibitem{May01}
Micaela Mayero.
\newblock {\em Formalisation et automatisation de preuves en analyses r\'eelle
  et num\'erique}.
\newblock PhD thesis, Universit\'e Paris VI, 2001.

\bibitem{ieee-754}
{Microprocessor Standards Committee}.
\newblock {IEEE} {S}tandard for {F}loating-{P}oint {A}rithmetic.
\newblock {\em IEEE Std. 754-2008}, pages 1--58, 2008.

\bibitem{or:vvs:10}
William~L. Oberkampf and Christopher~J. Roy.
\newblock {\em Verification and Validation in Scientific Computing}.
\newblock Cambridge University Press, 2010.

\bibitem{PetWilZeil96}
Marko Petkov\v{s}ek, Herbert~S. Wilf, and Doron Zeilberger.
\newblock {\em A=B}.
\newblock A K Peters Ltd, Wellesley, MA, 1996.
\newblock \url{http://www.cis.upenn.edu/~wilf/AeqB.html}.

\bibitem{str:fds:89}
John~C Strikwerda.
\newblock {\em Finite Difference Schemes and Partial Differential Equations}.
\newblock Chapman \& Hall, 1989.

\bibitem{tho:npd:95}
James~William Thomas.
\newblock {\em Numerical Partial Differential Equations: Finite Difference
  Methods}.
\newblock Number~22 in Texts in Applied Mathematics. Springer, 1995.

\bibitem{ver:rap:96}
R.~Verf\"urth.
\newblock {\em A Review of A Posteriori Error Estimation and Adaptive
  Mesh-Refinement Techniques}.
\newblock Teubner-Wiley, Stuttgart, Germany, 1996.

\bibitem{Zeil90}
Doron Zeilberger.
\newblock A fast algorithm for proving terminating hypergeometric identities.
\newblock {\em Discrete Math.}, 80:207--211, 1990.

\bibitem{Zeil91}
Doron Zeilberger.
\newblock The method of creative telescoping.
\newblock {\em J. Symbolic Computation}, 11:195--204, 1991.

\end{thebibliography}
\bibliographystyle{plain}

\appendix
\normalsize

\section{Glossary}
\label{sec:glossary}

This section gives the definition of concepts used in mathematical logic and
computer science that appear in this paper.

\begin{description}

\item[annotation]
  comment added to the {\clang}~code to specify the logical properties of the
  program.
  Tools turn them into \gloss{verification conditions}.

\item[atom]
  atomic component of a logical formula.
  In the setting of classical logic, any propositional formula can be rewritten
  using literals (atoms or negated atoms), conjunctions, and disjunctions
  only.

\item[automated theorem prover]
  software tool that automatically proves \gloss{goals}.
  It may fail to find a proof, even for valid formulas.

\item[decision procedure]
  an algorithm dedicated to proving specific properties.
  This is one of the basic blocks of theorem provers.

\item[deductive verification]
  process of verifying, with the help of theorem provers, that a program
  satisfies its \gloss{specification}.

\item[first-order logic]
  formal language of logical formulas that use quantifications over values
  only, and not over predicates and functions.

\item[formal proof]
  a finite sequence of deduction steps which are checked by a computer.

\item[forward error analysis]
  process of propagating error bounds from inputs to outputs of functions.

\item[functional programming]
  programming paradigm that treats computation as mathematical evaluation of
  functions and considers functions as ordinary values.

\item[goal]
  a set of hypotheses and a logical formula.
  The proof of a goal is a way to deduce the logical formula from the
  hypotheses.
  When proving a theorem, the statement of the theorem is the initial goal.

\item[higher-order logic]
  formal language of logical formulas that use quantifications over values,
  predicates and functions.

\item[inference rule]
  generic way of drawing a valid conclusion based on the form of hypotheses.
  Each deduction step of a \gloss{formal proof} must satisfy one of the
  inference rules of the logical system.

\item[interactive proof assistant]
  software tool to assist with the development of \gloss{formal proofs} by
  human-machine collaboration.

\item[interval arithmetic]
  arithmetic that operates on sets of values (typically intervals) instead of
  values.

\item[loop invariant]
  logical formula about the state of a program, that is valid before entering a
  loop and remains valid at the end of each iteration of the loop.

\item[semantics]
  meaning associated to each syntactic construct of a language.

\item[SMT solver]
  variety of \gloss{automated theorem prover} combining a SAT solver
  (propositional logic), equality reasoning, \gloss{decision procedures}
  ({\eg}, for linear arithmetic), and quantifier instantiation.

\item[specification]
  description of the expected behavior of a program.

\item[static analysis]
  analysis of a program without executing it.

\item[straight-line program]
  program that does not contain any loop and thus does not need any \gloss{loop
    invariant} to be verified.

\item[tactic]
  command for an \gloss{interactive proof assistant} to transform the current
  \gloss{goal} into one or more goals that imply it.

\item[validation]
  process of making sure of the correctness of a program by experiments such as
  tests.

\item[verification]
  process of making sure of the correctness of a program by mathematical means
  such as \gloss{formal proofs}.

\item[verification conditions]
  \gloss{goals} that need to be proved to guarantee the adequacy between the
  program and its \gloss{specification}.

\end{description}

\section{Fully Annotated Source Code}
\label{sec:annotated_source}
\begin{lstlisting}[style=ACSL, style=long]
/*@ axiomatic dirichlet_maths {
  @
  @ logic real c;
  @ logic real p0(real x);
  @ logic real psol(real x, real t);

  @ axiom c_pos: 0 < c;

  @ logic real psol_1(real x, real t);
  @ axiom psol_1_def:
  @   \forall real x; \forall real t;
  @   \forall real eps; \exists real C; 0 < C && \forall real dx;
  @   0 < eps ==> \abs(dx) < C ==>
  @   \abs((psol(x + dx, t) - psol(x, t)) / dx - psol_1(x, t)) < eps;

  @ logic real psol_11(real x, real t);
  @ axiom psol_11_def:
  @   \forall real x; \forall real t;
  @   \forall real eps; \exists real C; 0 < C && \forall real dx;
  @   0 < eps ==> \abs(dx) < C ==>
  @   \abs((psol_1(x + dx, t) - psol_1(x, t)) / dx - psol_11(x, t)) < eps;

  @ logic real psol_2(real x, real t);
  @ axiom psol_2_def:
  @   \forall real x; \forall real t;
  @   \forall real eps; \exists real C; 0 < C && \forall real dt;
  @   0 < eps ==> \abs(dt) < C ==>
  @   \abs((psol(x, t + dt) - psol(x, t)) / dt - psol_2(x, t)) < eps;

  @ logic real psol_22(real x, real t);
  @ axiom psol_22_def:
  @   \forall real x; \forall real t;
  @   \forall real eps; \exists real C; 0 < C && \forall real dt;
  @   0 < eps ==> \abs(dt) < C ==>
  @   \abs((psol_2(x, t + dt) - psol_2(x, t)) / dt - psol_22(x, t)) < eps;

  @ axiom wave_eq_0: \forall real x; 0 <= x <= 1 ==> psol(x, 0) == p0(x);
  @ axiom wave_eq_1: \forall real x; 0 <= x <= 1 ==> psol_2(x, 0) == 0;
  @ axiom wave_eq_2:
  @   \forall real x; \forall real t;
  @   0 <= x <= 1 ==> psol_22(x, t) - c * c * psol_11(x, t) == 0;
  @ axiom wave_eq_dirichlet_1: \forall real t; psol(0, t) == 0;
  @ axiom wave_eq_dirichlet_2: \forall real t; psol(1, t) == 0;

  @ logic real psol_Taylor_3(real x, real t, real dx, real dt);
  @ logic real psol_Taylor_4(real x, real t, real dx, real dt);

  @ logic real alpha_3; logic real C_3;
  @ logic real alpha_4; logic real C_4;

  @ axiom psol_suff_regular_3:
  @   0 < alpha_3 && 0 < C_3 &&
  @   \forall real x; \forall real t; \forall real dx; \forall real dt;
  @   0 <= x <= 1 ==> \sqrt(dx * dx + dt * dt) <= alpha_3 ==>
  @   \abs(psol(x + dx, t + dt) - psol_Taylor_3(x, t, dx, dt)) <=
  @     C_3 * \abs(\pow(\sqrt(dx * dx + dt * dt), 3));

  @ axiom psol_suff_regular_4:
  @   0 < alpha_4 && 0 < C_4 &&
  @   \forall real x; \forall real t; \forall real dx; \forall real dt;
  @   0 <= x <= 1 ==> \sqrt(dx * dx + dt * dt) <= alpha_4 ==>
  @   \abs(psol(x + dx, t + dt) - psol_Taylor_4(x, t, dx, dt)) <=
  @     C_4 * \abs(\pow(\sqrt(dx * dx + dt * dt), 4));

  @ axiom psol_le:
  @   \forall real x; \forall real t;
  @   0 <= x <= 1 ==> 0 <= t ==> \abs(psol(x, t)) <= 1;

  @ logic real T_max;
  @ axiom T_max_pos: 0 < T_max;

  @ logic real C_conv; logic real alpha_conv;
  @ lemma alpha_conv_pos: 0 < alpha_conv;
  @
  @ } */

/*@ axiomatic dirichlet_prog {
  @
  @  predicate analytic_error{L}
  @    (double **p, integer ni, integer i, integer k, double a, double dt)
  @    reads p[..][..];
  @
  @  lemma analytic_error_le{L}:
  @    \forall double **p; \forall integer ni; \forall integer i;
  @    \forall integer nk; \forall integer k;
  @    \forall double a; \forall double dt;
  @    0 < ni ==> 0 <= i <= ni ==> 0 <= k ==>
  @    0 < \exact(dt) ==>
  @    analytic_error(p, ni, i, k, a, dt) ==>
  @    \sqrt(1. / (ni * ni) + \exact(dt) * \exact(dt)) < alpha_conv ==>
  @    k <= nk ==> nk <= 7598581 ==> nk * \exact(dt) <= T_max ==>
  @    \exact(dt) * ni * c <= 1 - 0x1.p-50 ==>
  @    \forall integer i1; \forall integer k1;
  @    0 <= i1 <= ni ==> 0 <= k1 < k ==>
  @    \abs(p[i1][k1]) <= 2;
  @
  @  predicate separated_matrix{L}(double **p, integer leni) =
  @    \forall integer i; \forall integer j;
  @    0 <= i < leni ==> 0 <= j < leni ==> i != j ==>
  @    \base_addr(p[i]) != \base_addr(p[j]);
  @
  @ logic real sqr_norm_dx_conv_err{L}
  @   (double **p, real dx, real dt, integer ni, integer i, integer k)
  @   reads p[..][..];
  @ logic real sqr(real x) = x * x;
  @ lemma sqr_norm_dx_conv_err_0{L}:
  @   \forall double **p; \forall real dx; \forall real dt;
  @   \forall integer ni; \forall integer k;
  @   sqr_norm_dx_conv_err(p, dx, dt, ni, 0, k) == 0;
  @ lemma sqr_norm_dx_conv_err_succ{L}:
  @   \forall double **p; \forall real dx; \forall real dt;
  @   \forall integer ni; \forall integer i; \forall integer k;
  @   0 <= i ==>
  @   sqr_norm_dx_conv_err(p, dx, dt, ni, i + 1, k) ==
  @     sqr_norm_dx_conv_err(p, dx, dt, ni, i, k) +
  @     dx * sqr(psol(1. * i / ni, k * dt) - \exact(p[i][k]));
  @ logic real norm_dx_conv_err{L}
  @   (double **p, real dt, integer ni, integer k) =
  @   \sqrt(sqr_norm_dx_conv_err(p, 1. / ni, dt, ni, ni, k));
  @
  @ } */

/*@ requires leni >= 1 && lenj >= 1;
  @ ensures
  @   \valid_range(\result, 0, leni - 1) &&
  @   (\forall integer i; 0 <= i < leni ==>
  @    \valid_range(\result[i], 0, lenj - 1)) &&
  @   separated_matrix(\result, leni);
  @ */
double **array2d_alloc(int leni, int lenj);

/*@ requires (l != 0) && \round_error(x) <= 5./2*0x1.p-53;
  @ ensures
  @   \round_error(\result) <= 14 * 0x1.p-52 &&
  @   \exact(\result) == p0(\exact(x));
  @ */
double p_zero(double xs, double l, double x);

/*@ requires
  @   ni >= 2 && nk >= 2 && l != 0 &&
  @   dt > 0. && \exact(dt) > 0. &&
  @   \exact(v) == c && \exact(v) == v &&
  @   0x1.p-1000 <= \exact(dt) &&
  @   ni <= 2147483646 && nk <= 7598581 &&
  @   nk * \exact(dt) <= T_max &&
  @   \abs(\exact(dt) - dt) / dt <= 0x1.p-51 &&
  @   0x1.p-500 <= \exact(dt) * ni * c <= 1 - 0x1.p-50 &&
  @   \sqrt(1. / (ni * ni) + \exact(dt) * \exact(dt)) < alpha_conv;
  @
  @ ensures
  @   \forall integer i; \forall integer k;
  @   0 <= i <= ni ==> 0 <= k <= nk ==>
  @   \round_error(\result[i][k]) <= 78. / 2 * 0x1.p-52 * (k + 1) * (k + 2);
  @
  @ ensures
  @   \forall integer k; 0 <= k <= nk ==>
  @   norm_dx_conv_err(\result, \exact(dt), ni, k) <=
  @     C_conv * (1. / (ni * ni) + \exact(dt) * \exact(dt));
  @ */
double **forward_prop(int ni, int nk, double dt, double v,
                      double xs, double l) {

  /* Output variable. */
  double **p;

  /* Local variables. */
  int i, k;
  double a1, a, dp, dx;

  dx = 1./ni;
  /*@ assert
    @   dx > 0. && dx <= 0.5 &&
    @   \abs(\exact(dx) - dx) / dx <= 0x1.p-53;
    @ */

  /* Compute the constant coefficient of the stiffness matrix. */
  a1 = dt/dx*v;
  a = a1*a1;
  /*@ assert
    @   0 <= a <= 1 &&
    @   0 < \exact(a) <= 1 &&
    @   \round_error(a) <= 0x1.p-49;
    @ */

  /* Allocate space-time variable for the discrete solution. */
  p = array2d_alloc(ni+1, nk+1);

  /* First initial condition and boundary conditions. */
  /* Left boundary. */
  p[0][0] = 0.;
  /* Time iteration -1 = space loop. */
  /*@ loop invariant
    @   1 <= i <= ni &&
    @   analytic_error(p, ni, i - 1, 0, a, dt);
    @ loop variant ni - i; */
  for (i=1; i<ni; i++) {
    p[i][0] = p_zero(xs, l, i*dx);
  }
  /* Right boundary. */
  p[ni][0] = 0.;
  /*@ assert analytic_error(p, ni, ni, 0, a, dt); */

  /* Second initial condition (with p_one=0) and boundary conditions. */
  /* Left boundary. */
  p[0][1] = 0.;
  /* Time iteration 0 = space loop. */
  /*@ loop invariant
    @   1 <= i <= ni &&
    @   analytic_error(p, ni, i - 1, 1, a, dt);
    @ loop variant ni - i; */
  for (i=1; i<ni; i++) {
    /*@ assert \abs(p[i-1][0]) <= 2; */
    /*@ assert \abs(p[i][0]) <= 2; */
    /*@ assert \abs(p[i+1][0]) <= 2; */
    dp = p[i+1][0] - 2.*p[i][0] + p[i-1][0];
    p[i][1] = p[i][0] + 0.5*a*dp;
  }
  /* Right boundary. */
  p[ni][1] = 0.;
  /*@ assert analytic_error(p, ni, ni, 1, a, dt); */

  /* Evolution problem and boundary conditions. */
  /* Propagation = time loop. */
  /*@ loop invariant
    @   1 <= k <= nk &&
    @   analytic_error(p, ni, ni, k, a, dt);
    @ loop variant nk - k; */
  for (k=1; k<nk; k++) {
    /* Left boundary. */
    p[0][k+1] = 0.;
    /* Time iteration k = space loop. */
    /*@ loop invariant
      @   1 <= i <= ni &&
      @   analytic_error(p, ni, i - 1, k + 1, a, dt);
      @ loop variant ni - i; */
    for (i=1; i<ni; i++) {
      /*@ assert \abs(p[i-1][k]) <= 2; */
      /*@ assert \abs(p[i][k]) <= 2; */
      /*@ assert \abs(p[i+1][k]) <= 2; */
      /*@ assert \abs(p[i][k-1]) <= 2; */
      dp = p[i+1][k] - 2.*p[i][k] + p[i-1][k];
      p[i][k+1] = 2.*p[i][k] - p[i][k-1] + a*dp;
    }
    /* Right boundary. */
    p[ni][k+1] = 0.;
    /*@ assert analytic_error(p, ni, ni, k + 1, a, dt); */
  }

  return p;

}
\end{lstlisting}

\section{Fundamental Solution and Jacobi Polynomials}
\label{sec:jacobi}

Let~$\lambda$ be the sequence defined in
Section~\ref{sec:global-round-off-error}.
Note that adding zero initial values for the fictitious time step $k=-1$
makes this sequence be a time shift by one of the fundamental solution of the
discrete acoustic wave equation, associated with the input
data~$s_{\rm h}\equiv 0$, $u_{0,{\rm h}}\equiv 0$, for all~$i$, $i\not=0$,
$u_{1,i}=0$, and $u_{1,0}=1$.
The items of the sequence satisfy the following equations:
\begin{eqnarray}
  \label{eq:init0}
  \forall i, & & \lambda_i^{-1} = 0, \\
  \label{eq:init1}
  \forall i \not= 0, & & \lambda_i^0 = 0, \quad \lambda_0^0 = 1, \\
  \label{eq:iter}
  \forall i, \forall k \geq 0, & &
  \lambda_i^{k+1} =
  a (\lambda_{i-1}^k + \lambda_{i+1}^k)
  + 2 (1 - a) \lambda_i^k - \lambda_i^{k-1}.
\end{eqnarray}
We want to prove Lemma~\ref{th:jacobi}, {\ie} that for all $i,k$,
$k\geq 0$, we have $\lambda_i^k\geq 0$.

The proof is highly indebted to computer algebra: we use a generating
function to obtain a closed-form expression for the $\lambda$'s, the creative
telescoping method of Zeilberger~\cite{PetWilZeil96} to express those
$\lambda$'s in terms of Jacobi polynomials, and finally a result by Askey and
Gasper~\cite{AsGa72} to ensure the nonnegativity.
We have not mechanically checked this proof.
For example, the Askey and Gasper result would have required enormous {\coq}
developments, but parts of it could have been formalized, in particular
Zeilberger's algorithm provides a certificate that eases the
\gloss{verification} of its result.

Consider the associated bivariate generating function formally defined by
\begin{equation*}
  \Lambda(X, T) = \sum_i \sum_{k \geq -1} \lambda_i^k X^i T^k.
\end{equation*}
The above recurrence relation in Equation~(\ref{eq:iter}) rewrites
\begin{eqnarray*}
  \sum_i \sum_{k \geq 0} \lambda_i^{k+1} X^i T^k & = &
  a \left( \sum_i \sum_{k \geq 0} \lambda_{i-1}^k X^i T^k
    + \sum_i \sum_{k \geq 0} \lambda_{i+1}^k X^i T^k \right) \\
  & & + 2 (1 - a) \sum_i \sum_{k \geq 0} \lambda_i^k X^i T^k
  - \sum_i \sum_{k \geq 0} \lambda_i^{k-1} X^i T^k.
\end{eqnarray*}
Since coefficients for $k$'s smaller than~1 are almost all equal to zero,
we formally have
\begin{equation*}
  \begin{array}{ll}
    \displaystyle
    T \sum_i \sum_{k \geq 0} \lambda_i^{k+1} X^i T^k = \Lambda(X, T) - 1, &
    \displaystyle
    \sum_i \sum_{k \geq 0} \lambda_{i-1}^k X^i T^k = X \Lambda(X, T), \\
    \displaystyle
    X \sum_i \sum_{k \geq 0} \lambda_{i+1}^k X^i T^k = \Lambda(X, T), &
    \displaystyle
    \sum_i \sum_{k \geq 0} \lambda_i^k X^i T^k = \Lambda(X, T), \\
    \multicolumn{2}{c}{%
      \displaystyle
      \sum_i \sum_{k \geq 0} \lambda_i^{k-1} X^i T^k = T \Lambda(X, T).}
  \end{array}
\end{equation*}
Therefore, the generating function satisfies
\begin{eqnarray*}
  \frac{\Lambda(X, T) - 1}{T} & = &
  a \left(X + \frac{1}{X}\right) \Lambda(X, T)
  + 2 (1 - a) \Lambda(X, T) - T \Lambda(X, T),
\end{eqnarray*}
which solution is given by
\begin{eqnarray*}
  \Lambda(X, T) & = &
  \frac{1}{\displaystyle
    (1 - T)^2 \left[
      1 - a \frac{T}{X} \left(\frac{1 - X}{1 - T}\right)^2
    \right]}.
\end{eqnarray*}

Now, we can evaluate the power series expansion of the generating
function~$\Lambda$.
Using the following power series expansion, valid for $|u|<1$,
\begin{eqnarray*}
  \frac{1}{(1 - u)^{p+1}} & = & \sum_{n\geq 0} \binom{p+n}{p} u^n,
\end{eqnarray*}
and the properties of binomial coefficients, we successively have
\begin{eqnarray*}
  \Lambda(X, T) & = &
  \frac{1}{(1 - T)^2}
  \sum_{n\geq 0} a^n \frac{T^n}{X^n} \left(\frac{1 - X}{1 - T}\right)^{2n} \\
  & = &
  \sum_{n \geq 0} a^n
  \sum_{i=0}^{i=2n} \binom{2n}{i} (-1)^i X^{i-n}
  \sum_{k\geq 0} \binom{2n+1+k}{2n+1} T^{n+k} \\
  & = &
  \sum_i X^i \sum_{k \geq 0} T^k
  \sum_{n=|i|}^{n=k} \binom{2n}{n+i} \binom{n+k+1}{2n+1} (-1)^{n+i} a^n.
\end{eqnarray*}
Finally, by identification of the two power series expansions of the generating
function~$\Lambda$, we have for all $i,k$, $0\leq|i|\leq k$,
\begin{equation}
  \label{eq:lambda_i_k}
  \lambda_i^k =
  \sum_{n=|i|}^{n=k} \binom{2n}{n+i} \binom{n+k+1}{2n+1} (-1)^{n+i} a^n.
\end{equation}

For $i=0$, sharp eyes may recognize that $\lambda_0^k=\sum_{n=0}^kP_n(1-2a)$
where the~$P_n$'s are Legendre polynomials.
Fej\'er showed in~\cite{Fejer08} the nonnegativity of the sum of Legendre
polynomials when the argument is in $[-1,1]$, which is satisfied here since
we consider $0<a<1$.
More generally, we have $\lambda_i^k=a^{|i|}\sum_{n=0}^{k-|i|}P_n^{(2|i|,0)}(1-2a)$
where the~$P_n^{(\alpha,\beta)}$'s are Jacobi polynomials.
Askey and Gasper generalized in~\cite{AsGa72,Gasper77} Fej\'er's result for
$\beta\geq 0$ and $\alpha+\beta\geq -2$.
See also~\cite{AnAsRo99}, pages~314 and~384.

Indeed, from the definition of Jacobi polynomials, for all $\alpha,\beta>-1$,
for all $n\in\N$, for all $x\in[-1,1]$,
\begin{equation*}
  P_n^{(\alpha, \beta)} (x) =
  \sum_{p=0}^n \binom{n + \alpha}{p} \binom{n + \beta}{n - p}
  \left( \frac{x + 1}{2} \right)^p
  \left( \frac{x - 1}{2} \right)^{n - p},
\end{equation*}
we have, for all $i,k$, $0\leq i\leq k$,
\begin{eqnarray*}
  a^i\sum_{n=0}^{k-i}P_n^{(2i,0)}(1-2a) & = &
  a^i \sum_{n=0}^{k-i} \sum_{p=0}^n \binom{n + 2i}{p} \binom{n}{n - p}
  (1 - a)^p (-a)^{n-p} \\
  & = &
  \sum_{n=i}^k \sum_{p=0}^{n-i} \sum_{q=0}^p
  \binom{n + i}{p} \binom{n - i}{p} \binom{p}{q}
  (-1)^{n-p+q+i} a^{n-p+q} \\
  & = &
  \sum_{n=i}^k \sum_{p=i}^n \sum_{q=p}^n
  \binom{n + i}{n - p} \binom{n - i}{n - p} \binom{n - p}{n - q}
  (-1)^{q+i} a^q \\
  & = &
  \sum_{n=i}^k \sum_{p=i}^n \sum_{q=n}^k
  \binom{q + i}{p + i} \binom{q - i}{p - i} \binom{q - p}{n - p}
  (-1)^{n+i} a^n.
\end{eqnarray*}
We have successively shifted~$n$ by~$i$, replaced $n-p$ by~$p$, then
shifted~$q$ by~$p$.
To obtain the last equality, notice that the previous triple sum is actually
taken over $\{(n,p,q)\in\N^3/i\leq p\leq q\leq n\leq k\}$, hence we can
take~$q$ in $[i..k]$, $p$ in $[i..q]$, $n$ in $[q..k]$, and then switch
notations~$n$ and~$q$, and use the symmetry of binomial coefficients.
Identifying with the expression of Equation~(\ref{eq:lambda_i_k}), we are led
to prove, for all $i,n,k$, $0\leq i\leq n\leq k$,
\begin{equation}
  \sum_{p=i}^n \sum_{q=n}^k
  \binom{q + i}{p + i} \binom{q - i}{p - i} \binom{q - p}{n - p} =
  \binom{2n}{n+i} \binom{n+k+1}{2n+1}.
\end{equation}

Suppose we have the following identity, for all $i,n,k$, $0\leq i\leq n\leq k$,
\begin{equation}
  \label{eq:hard_one}
  \sum_{p=i}^n
  \binom{k + i}{p + i} \binom{k - i}{p - i} \binom{k - p}{n - p} =
  \binom{2n}{n+i} \binom{k+n}{2n}.
\end{equation}
Then, we would have
\begin{multline*}
  \sum_{p=i}^n \sum_{q=n}^k
  \binom{q + i}{p + i} \binom{q - i}{p - i} \binom{q - p}{n - p} = \\
  \sum_{q=n}^k \binom{2n}{n+i} \binom{q+n}{2n}
  = \binom{2n}{n+i} \sum_{q=2n}^{k+n} \binom{q}{2n}
  = \binom{2n}{n+i} \binom{k+n+1}{2n+1}.
\end{multline*}
The last equality comes directly from the recurrence relation for the binomial
coefficients (column-sum property of Pascal's triangle).

Proving identity of Equation~(\ref{eq:hard_one}) is a bit more technical.
The hypergeometric nature of its terms makes it a good candidate for
Zeilberger's algorithm, a.k.a. the method of creative telescoping,
see~\cite{Zeil90,Zeil91}.
Let us introduce some new notations, for all $i,n,k$, $0\leq i\leq n\leq k$,
\begin{eqnarray*}
  F (i, n, k; p) & = &
  \binom{k + i}{p + i} \binom{k - i}{p - i} \binom{k - p}{n - p}, \\
  f (i, n, k) & = & \sum_{p} F (i, n, k; p), \\
  g (i, n, k) & = & \binom{2n}{n+i} \binom{k+n}{2n}.
\end{eqnarray*}
Note that~$F$ vanishes when~$p$ is outside the interval $[i..n]$.
Thus, identity of Equation~(\ref{eq:hard_one}) now writes $f=g$.

Let~$I$ (resp. $N$, $K$, and~$P$) be the forward shift operator in~$i$
(resp. in~$n$, $k$, and~$p$).
E.g., $If(i,n,k)=f(i+1,n,k)$.
We first assume that the function~$f$ satisfies the following first-order
recurrence relations, for all $i,n,k$, $0\leq i\leq n\leq k$,
\begin{eqnarray}
  \label{eq:recf_i}
  (n+1+i) I f & = & (n-i) f, \\
  \label{eq:recf_n}
  (n+1+i) (n+1-i) N f & = & (k+1+n) (k-n) f, \\
  \label{eq:recf_k}
  (k+1-n) K f & = & (k+1+n) f.
\end{eqnarray}
Then, it is easy to show that the function~$g$ satisfies exactly the same
first-order recurrence relations, and since $f(0,0,0)=g(0,0,0)=1$, we have
$f=g$.
Indeed, simply using the symmetry and absorption properties of binomial
coefficients, we have, for all $i,n,k$, $0\leq i\leq n\leq k$,
\begin{equation*}
  \frac{g (i+1, n, k)}{g (i, n, k)}
  = \frac{\binom{2n}{n+i+1}}{\binom{2n}{n+i}}
  = \frac{n-i}{n+1+i},
\end{equation*}
\begin{multline*}
  \frac{g (i, n+1, k)}{g (i, n, k)} =
  \mymultlinenewline
  \frac{\binom{2n+2}{n+1+i} \binom{k+n+1}{2n+2}}%
         {\binom{2n}{n+i} \binom{k+n}{2n}}
  = \frac{k+1+n}{n+1+i} \frac{\binom{2n+1}{n+i}}{\binom{2n}{n+i}}
    \frac{\binom{k+n}{2n+1}}{\binom{k+n}{2n}}
  = \frac{(k+1+n)(k-n)}{(n+1+i)(n+1-i)},
\end{multline*}
\begin{equation*}
  \frac{g (i, n, k+1)}{g (i, n, k)}
  = \frac{\binom{k+1+n}{2n}}{\binom{k+n}{2n}}
  = \frac{k+1+n}{k+1-n}.
\end{equation*}

Finding the first-order recurrence relations for~$f$,
{\ie} Equations~(\ref{eq:recf_i}), (\ref{eq:recf_n}), and~(\ref{eq:recf_k}), is
the job of the method of creative telescoping.
Actually, Zeilberger's algorithm provides recurrence relations for the
hypergeometric summand~$F$ of the form
\begin{equation*}
  \sum_{m=0}^{m^\prime} b_{l,m} L^m F = (P - 1) (R_l F)
\end{equation*}
where coefficients~$b_{l,m}$ are polynomials independent of~$p$, and~$R_l$ is a
rational function.
There are actually one such recurrence relation per free variable~$l$
(here~$i$, $n$, and~$k$), and $L$ is the generic forward shift operator in the
generic free variable~$l$.
Thus, since coefficients~$b_{l,m}$ do not depend on~$p$, when summing over~$p$,
the right-hand terms telescope, and we can deduce similar recurrence relations
for the sum~$f$
\begin{equation*}
  \sum_{m=0}^{m^\prime} b_{l,m} L^m f = 0.
\end{equation*}
Note that although those recurrence relations are difficult to obtain, and even
to check on the sum~$f$, they are easy to check on the summand~$F$ (since there
is no more sum over~$p$).
On has just to check the simpler equations with rational expressions
\begin{equation}
  \label{eq:recF}
  b_{l,0} + \sum_{m=1}^{m^\prime} b_{l,m} \frac{L^m F}{F} =
  P R_l \frac{P F}{F} - R_l.
\end{equation}
In the present case, Zeilberger's algorithm provides first-order recurrence
relations with
\begin{equation*}
  \begin{array}{lll}
    b_{i,0} = n-i, & b_{n,0} = (k+1+n)(k-n), & b_{k,0} = k+1+n \\
    b_{i,1} = -(n+1+i), & b_{n,1} = -(n+1+i)(n+1-i), & b_{k,1} = -(k+1-n), \\
    \displaystyle
    R_i = \frac{(1 + 2i)(p - i)}{k - i}, &
    \displaystyle
    R_n = \frac{(k - n)(p + i)(p - i)}{n + 1 - p}, &
    \displaystyle
    R_k = \frac{(p + i)(p - i)}{k + 1 - p}.
  \end{array}
\end{equation*}
And, simply using the symmetry and absorption properties of binomial
coefficients, Equation~(\ref{eq:recF}) is successively equivalent to,
for $l=i$,
\begin{eqnarray*}
  (\ref{eq:recF})
  & \Leftrightarrow &
  (n-i) - (n+1+i) \frac{(k+1+i)}{(p+1+i)} \frac{(p-i)}{(k-i)} \\
  & &
  = \frac{(1+2i)(p+1-i)}{(k-i)}
  \frac{(k-p)}{(p+1+i)} \frac{(k-p)}{(p+1-i)} \frac{(n-p)}{(k-p)}
  \myeqnarraynewline
  - \frac{(1+2i)(p-i)}{(k-i)} \\
  & \Leftrightarrow &
  (n-i)(k-i)(p+1+i) - (n+1+i)(k+1+i)(p-i) \\
  & &
  = (1+2i)(k-p)(n-p) - (1+2i)(p-i)(p+1+i) \\
  & \Leftrightarrow &
  \frac{(1+2i)}{4} \left[
    (2n+1)(2k+1) - (2n+1)(2p+1) - (2k+1)(2p+1) \right. \\
  & & \left. + (1+2i)^2 \right]
  = (1+2i) \left[ (k-p)(n-p) - (p-i)(p+1+i) \right] \\
  & \Leftrightarrow &
  \frac{(1+2i)^2}{4}- \frac{(2p+1)^2}{4} = - (p-i)(p+1+i),
\end{eqnarray*}
for $l=n$,
\begin{eqnarray*}
  (\ref{eq:recF})
  & \Leftrightarrow &
  (k+1+n)(k-n) - (n+1+i)(n+1-i) \frac{(k-n)}{(n+1-p)} \\
  & &
  = \frac{(k-n)(p+1+i)(p+1-i)}{(n-p)}
  \frac{(k-p)}{(p+1+i)} \frac{(k-p)}{(p+1-i)} \frac{(n-p)}{(k-p)} \\
  & & - \frac{(k-n)(p+i)(p-i)}{(n+1-p)} \\
  & \Leftrightarrow &
  (n+1+k)(n+1-p) - (n+1+i)(n+1-i) \\
  & & = (k-p)(n+1-p) - (p+i)(p-i) \\
  & \Leftrightarrow &
  (n+1)(k-p) - kp + i^2 = (k-p)(n+1) -kp + i^2,
\end{eqnarray*}
and for $l=k$,
\begin{eqnarray*}
  (\ref{eq:recF})
  & \Leftrightarrow &
  (k+1+n) - (k+1-n)
  \frac{(k+1+i)}{(k+1-p)} \frac{(k+1-i)}{(k+1-p)} \frac{(k+1-p)}{(k+1-n)} \\
  & &
  = \frac{(p+1+i)(p+1-i)}{(k-p)}
  \frac{(k-p)}{(p+1+i)} \frac{(k-p)}{(p+1-i)} \frac{(n-p)}{(k-p)}
  \myeqnarraynewline
  - \frac{(p+i)(p-i)}{(k+1-p)} \\
  & \Leftrightarrow &
  (k+1+n)(k+1-p) - (k+1+i)(k+1-i) \\
  & & = (n-p)(k+1-p) - (p+i)(p-i) \\
  & \Leftrightarrow &
  (k+1)(n-p) - np + i^2 = (n-p)(k+1) -np + i^2.
\end{eqnarray*}
Which are all true.

\end{document}